# An Exact Convex Formulation of the Optimal Power Flow in Radial Distribution Networks Including Transverse Components

Mostafa Nick, Rachid. Cherkaoui, Jean-Yves Le Boudec, and Mario Paolone
École Polytechnique Fédérale de Lausanne (EPFL)

***Abstract*—The recent literature discusses the use of the relaxed Second Order Cone Programming (SOCP) for formulating Optimal Power Flow problems (OPF) for radial power grids. However, if the shunt parameters of the lines that compose the power grid are considered, the proposed methods do not provide sufficient conditions that can be verified ex ante for the exactness of the optimal solutions. The same formulations also have not correctly accounted for the lines' ampacity constraint. Similar to the inclusion of upper voltage-magnitude limit, the SOCP relaxation faces difficulties when the ampacity constraints of the lines are binding. In order to overcome these limitations, we propose a convex formulation for the OPF in radial power grids, for which the AC-OPF equations, including the transverse parameters, are considered. To limit the lines' current together with the nodal voltage-magnitudes, we augment the formulation with a new set of more conservative constraints. Sufficient conditions are provided to ensure the feasibility and optimality of the proposed OPF solution. Furthermore, the proofs of the exactness of the SOCP relaxation are provided. Using standard power grids, we show that these conditions are mild and hold for real distribution networks.**

## I. INTRODUCTION

THE Optimal Power Flow (OPF) is a well-known challenging optimization problem. It is the main building block for the formulation of optimal controls, as well as operation and planning problems in power systems. Typical examples are unit commitment, grid planning and reactive-power dispatch problems [1]-[4].

The OPF is inherently a non-convex optimization problem, consequently its solution is challenging (e.g., [5]). The authors in [6] show that the AC-feasibility problem (finding a feasible solution for OPF problem) is NP-hard for radial networks. Several methods have been proposed to solve OPF. We classify them into categories: (i) approximated methods that modify the physical description of the power flows equations [7]-[9], (ii) non-linear optimization methods [10]-[12], (iii) heuristic methods [13], [14], and (iv) convexification approaches [15]-[29]. Here, we briefly recall the main characteristics of the convex models of the AC-OPF, as they are the most representative. Several relaxations have been applied to convexify the AC-OPF problem (see the survey discussed in [15]). The authors of [16] propose a linearized-relaxed model to approximate and solve the AC-OPF problem. A semidefinite relaxation method is proposed in [17]. The authors in [18] and [19] investigate the application of moment-based relaxation for the OPF solution. Several recent papers have focused on the OPF problem in radial unbalanced power grids [20]-[23]. In [20] and [21] it is proposed to use SDP relaxation for the solution of this non-convex optimization problem. A distributed optimization, based on the alternative direction method of multipliers (ADMM) and semidefinite relaxation, is proposed in [22]. In [23] the authors propose a dedicated distributed optimization procedure, still based on the ADMM, for the OPF solution in unbalanced radial grids. In order to decrease the computational complexity, they reduce the ADMM subproblems to either a closed form solution or an Eigen-decomposition of a $6 \times 6$ Hermitian matrix. Even if these papers have proposed technically deployable techniques, they did not formally prove the exactness of their relaxations, as only numerical verifications have been provided for specific grids.

In case of balanced radial grids, the Second Order Cone Programming (SOCP) relaxation is proposed in [24]-[28] to solve the OPF.

In this paper, we address this last category of OPF solution methods (OPF in tree networks). We refer to this specific category because of the increasing needs of Distribution Network Operators (DNOs) to actively control their grids in an optimal fashion due to the increasing connection of the distributed generation (mainly from renewable energy resources) and the controllable devices such as distributed storage and demand response.

The authors in [25] investigate the geometry of the feasible injection region in radial distribution networks. In particular, they show that the SOCP (also SDP) relaxation is exact when (i) the angle separation of voltage phasors of line terminals is sufficiently small, and (ii) there is no bound on nodal reactive power (this last is a condition that is hard to meet for real systems). The latest contribution published on this subject is [28]. The authors show that the SOCP relaxation is computationally more efficient than the semidefinite relaxation. In [26] the authors show that the SOCP relaxation is tight when there is no upper bound on the nodal load consumptions. In [27] it is shown that the relaxation is exact with no upper bound on the nodal voltage-magnitudes, together with a specific condition involving the network parameters that can be checked a priori. In [28] the authors improve the work of [27] and introduce a more conservative constraint on the upper bound for the nodal voltage-magnitudes.

Although the model proposed in [28] works properly in many operating conditions, it has two important shortcomings: it does not take into account the shunt capacitors of the equivalent two-port Π line model and the ampacity constraint of the lines (which is an important limitation, for instance, in grids with coaxial underground cables) [29]. All the proofs in [28] are provided without considering these two important elements. Regarding the lines' ampacity limits, as shown in [29] and in Section IV-A, this is a fundamental problem of the relaxations proposed in [28]. They cannot be addressed by simply adding more details to the model; essentially the relaxation becomes inexact when ampacity limits of the lines are taken into account.

The authors in [30] propose a methodology for addressing the above shortcoming. It is based on the augmented Lagrangian method for solving the original non-convex OPF problem that considers the shunt impedances and the ampacity limit of the lines. However, this method is iterative and, as a consequence, computationally expensive. A positive aspect is that it can be easily formulated in a distributed manner.

In this paper, our contributions are (i) to modify the OPF problem by adding a new set of constraints to ensure the feasibility of its optimal solution, (ii) to provide sufficient conditions, by taking into account the shunt elements of the lines, to ensure the exactness of the SOCP relaxation that could be verified ex-ante, and (iii) to prove the exactness of the proposed relaxed OPF model under the provided sufficient conditions. In other words, suppose $\mathcal{F}$ is the feasible set of power injections of OPF. We first create an inner approximation of $\mathcal{F}$, called $\mathcal{F}^{in}$. Then, we take an outer relaxation of $\mathcal{F}^{in}$, say $\mathcal{F}^{in_r}$. Finally, we show that $\mathcal{F}^{in_r}$ is exact with respect to $\mathcal{F}^{in}$ under mild conditions. Thus, we are able to recover a feasible point of the original OPF problem. In this respect, a fundamental difference with respect to the existing literature is that we obtain for any feasible solution of the relaxed optimization problem --not only the optimal one -- the set of power injections that corresponds to one solution of the original non-relaxed OPF problem.

We show that the proposed formulation is characterized by a slightly reduced space of feasible solutions. The removed space is normally close to the technical limits of the grid; this space is not a desirable operating region for DNOs. We also show that the sufficient conditions always hold when the line impedances are not too large, which is expected in most distribution networks.

The structure of the paper is the following: In Section II, we introduce the notation and the OPF problem. Then, we introduce the auxiliary variables and constraints in Section III. In Section IV, we illustrate the proposed formulations for OPF in radial networks. In Section V, we introduce the power flow equations in matrix form. We provide our theorem, as well as the exactness conditions and proof of exactness, in Section VI. In Section VII, to quantify the advantages of the proposed formulation and the margins where the provided conditions hold, we use the IEEE 34-bus test-case network and the CIGRE European benchmark medium-voltage network. Furthermore, in Section VIII, we demonstrate the capabilities of the proposed formulation to provide feasible and optimal solutions. We provide a discussion on the extension of the work to unbalanced radial grids. In Section IX, we conclude the paper with the final remarks, summarizing the main findings of this manuscript.

## II. Proposed Modified OPF Formulation

### A. Notations and Definitions

The network is radial. Index 0 is for the slack bus and its voltage is fixed ($v_0$). Without loss of generality, we can assume that only bus 1 is connected to the slack bus (otherwise the problem separates into several independent problems). Buses other than the slack bus are denoted with $1, \dots, L$; $\mathcal{L}$ denotes the set $\mathcal{L} = \{1,2,\dots,L\}$ and $\mathrm{up}(l)$ is the label of the bus that is upstream of bus $l$. We also label with $l$ the line whose downstream bus is bus $l$; its upstream bus is therefore $\mathrm{up}(l)$. For a set of lines $\mathcal{M} = \{1,2,\dots,M\}$, we denote the lines for which $\mathbf{H}_{l,m} = 1, \forall\, m \in \mathcal{M}$ upstream lines/buses, $\left(\mathcal{L}^{\mathcal{M}}\right)$ ($\mathbf{H}$ is defined in (1)). Finally, $\underline{\mathcal{L}}$ denotes the set of the buses that are the leaves of the graph that maps the grid topology.

Let $S_l^t = P_l^t + jQ_l^t$ be the complex power flow entering line $l$ from the top, i.e., from bus $\mathrm{up}(l)$; $S_l^c = P_l^t + jQ_l^c$ the complex power flow entering the central element of line $l$ (it is equal to $S_l^t$ minus the reactive power associated with the shunt admittance connected to bus $\mathrm{up}(l)$, see Fig. 1); $S_l^b = P_l^b + jQ_l^b$ the complex power flow entering bus $l$ from the bottom part of line $l$; and $f_l$ the square of the current in the central element of line $l$ (Fig.1). Let $z_l = r_l + jx_l$ and $2b_l$ be the longitudinal impedance and shunt capacitance of line $l$. We denote with $z_l^*$ the complex conjugate of $z_l$. We assume that $r_l, x_l$ and $b_l$ of the all lines are positive.

Let $v_l$ be the square of voltage magnitude and $s_l = p_l + jq_l$ be the power absorption at bus $l$. $p_l \geq 0$ and $q_l \geq 0$ denote power consumptions, $p_l \leq 0$ and $q_l \leq 0$ denote powers injections. Let $2B_l$ denote the sum of the susceptances of the lines connected to bus $l$.

$v^{\max}$ and $v^{\min}$ are the square of maximum and minimum magnitude of nodal voltages. $I_l^{\max}$ is the square of maximum current limit of line $l$ ($l \in \mathcal{L}$). $\Re(.)$ and $\Im(.)$ denote the real and imaginary parts of complex numbers, and $j \coloneqq \sqrt{-1}$ is the imaginary unit; $\max\{a, b\}$ returns the maximum of $a$ and $b$ and $\min\{a, b\}$ returns the minimum of $a$ and $b$.

A notation without subscript, such as $v$, denotes a column vector with $L$ rows as

$$v = \begin{pmatrix} v_1 \\ \vdots \\ v_L \end{pmatrix}, S = \begin{pmatrix} S_1^t \\ \vdots \\ S_L^t \end{pmatrix}, P = \begin{pmatrix} P_1^t \\ \vdots \\ P_L^t \end{pmatrix}, I^{\max} = \begin{pmatrix} I_1^{\max} \\ \vdots \\ I_L^{\max} \end{pmatrix}, \text{etc.}$$

Note that for $S, P, Q$ and their related auxiliary variables ($\bar{S}, \hat{S}$, ...) the vectors $S, P, Q$ represent the relevant values at the upper side of line $l$ ($S_l^t$). The notation $|P|$ represents the column vector with $L$ rows whose $l^{\text{th}}$ element is the absolute value $|P_l|$. The comparison of vectors is entry-wise, i.e., $P \leq P'$ means $P_l \leq P_l'$ for every $l \in \mathcal{L}$. The transposed of $P$ is denoted with $P^T$.

Matrices are shown with bold non-italic capital letters such as $\mathbf{A}$. We use the Euclidean (Frobenius) norm for vectors ($\|v\| = \sqrt[2]{\sum_{k=1}^{L}(v_k^2)}$) and also the Frobenius norm $\|\mathbf{A}\|$ for matrices (($\|\mathbf{A}\| = \sqrt[2]{\sum_{i=1}^{L}\left(\sum_{j=1}^{L}\left(a_{ij}^2\right)\right)}$ ). For two matrices $\mathbf{A}, \mathbf{B}$ of equal dimensions, the notation $\mathbf{A} \circ \mathbf{B}$ denotes their Hadamard product, defined by $(\mathbf{A} \circ \mathbf{B})_{k,l} = \mathbf{A}_{k,l}\mathbf{B}_{k,l}$ for all $k, l$.

For the reader's convenience, the matrices defined in the paper are listed below.

- $\mathbf{I}$ is the $L \times L$ identity matrix.
- For a vector such as $r$, $\operatorname{diag}(r)$ denotes the diagonal matrix whose $l^{\text{th}}$ element is $r_l$.
- $\mathbf{G}$ is the adjacency matrix of the oriented graph of the network, i.e. $\mathbf{G}_{k,l}$ is defined for $k, l \in \mathcal{L}$ and $\mathbf{G}_{k,l} = 1$ if $k = \operatorname{up}(l)$ and 0 otherwise. Diagonal elements are zero.
- $\mathbf{H}$ is the closure of $\mathbf{G}$, i.e. $\mathbf{H}_{k,l} = 1$ if bus $k$ is on the path from the slack bus to bus $l$ or $k = l$, and $\mathbf{H}_{k,l} = 0$ otherwise. Because the network is radial, $\mathbf{G}^L = 0$ and

$$\mathbf{H} = \mathbf{I} + \mathbf{G} + \mathbf{G}^2 + \cdots + \mathbf{G}^{L-1} = (\mathbf{I} - \mathbf{G})^{-1} \tag{1}$$

- $\mathbf{M} = \mathbf{2}\operatorname{diag}(x)\, \mathbf{H}\, \operatorname{diag}(B)$.
- $\mathbf{C} = (\mathbf{I} - \mathbf{G}^T - \mathbf{M})^{-\mathbf{1}}$. ($\mathbf{C}$ is well-defined and is non-negative (entry-wise) when condition *C1* (later defined in Section VI.A) holds).

- $\mathbf{D}$ is the entry-wise non-negative matrix defined by

$$\mathbf{D} = \mathbf{C}\big[2\mathrm{diag}(r)\big((\mathbf{H}-\mathbf{I})\mathrm{diag}(r)\big) + 2\mathrm{diag}(x)\big((\mathbf{H}-\mathbf{I})\mathrm{diag}(x)\big) + \mathrm{diag}(|z|^2)\big] \tag{2}$$

- $\pi, \varrho$ and $\vartheta$ are the vectors defined by

$$\pi_l = \frac{\max\left\{P_l^{\max}, \left|\mathbf{H}p^{\min}\right|_l\right\}}{v_l^{\min}} \tag{3}$$

$$\varrho_l = \frac{\max\left\{Q_l^{\max} + b_l v^{\max}, \left|\mathbf{H}q^{\min} - \mathbf{H}\mathrm{diag}(b)(\mathbf{I}+\mathbf{G}^{\mathbf{T}})(v^{\max})\right|_l\right\}}{v_l^{min}} \tag{4}$$

$$\vartheta_l = (\pi_l)^2 + (\varrho_l)^2 \tag{5}$$

where $p^{\min}$ and $q^{\min}$are the vectors of minimum absorptions level on the buses of the system. $P_l^{\max}$ and $Q_l^{\max}$ are the maximum allowed active and reactive power-flows of line $l$.

- $\mathbf{E}$ and $\mathbf{F}$ are the entry-wise non-negative matrices defined by:

$$\mathbf{F} = (\mathbf{H}\mathrm{diag}(x) + \mathbf{H}\mathrm{diag}(B)\mathbf{D}) \tag{6}$$

$$\mathbf{E} = 2\mathrm{diag}(\pi)\mathbf{H}\mathrm{diag}(r) + 2\mathrm{diag}(\varrho)\mathbf{F} + \mathrm{diag}(\vartheta)\mathbf{D} \tag{7.}$$

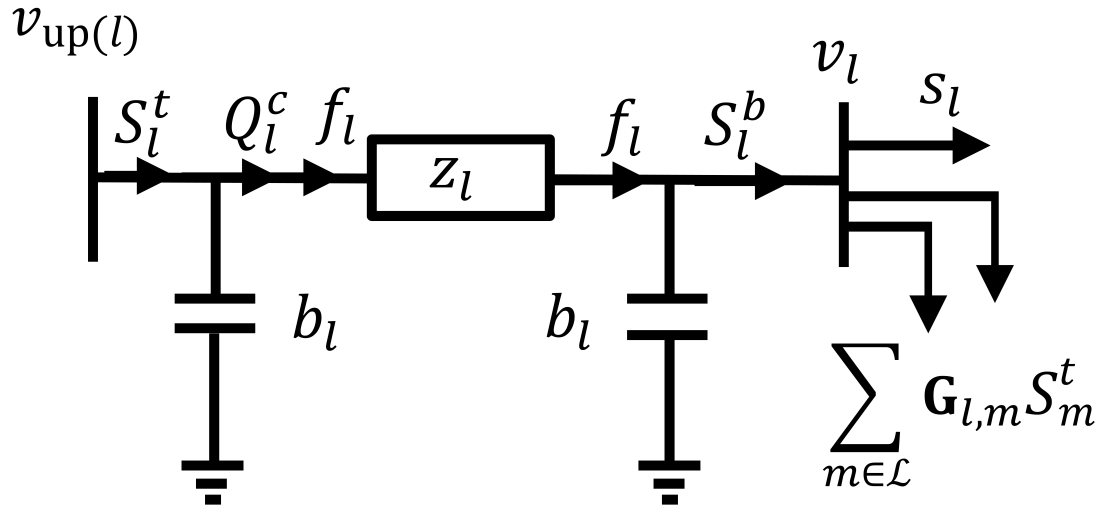

Fig.1. classic two-port Π model of a transmission line adopted for the formulation of the OPF relaxed constraints.

### *B. Power-Flow Equations with Transverse Parameters*

In this subsection, we introduce the power-flow equations inferred from the transmission line two-port Π model written with the notation of subsection II.A. For the sake of clarity, the transmission line two-port Π model is shown in Fig. 1.

For a given radial power network, the power-flow equations are given by (8).

$$S_l^t = s_l + \sum_{m\in\mathcal{L}} \mathbf{G}_{l,m} S_m^t + z_l f_l - j\big(v_{\mathrm{up}(l)} + v_l\big)b_l,\ \ \forall\, l \in \mathcal{L} \tag{8.a}$$

$$v_l = v_{\mathrm{up}(l)} - 2\Re\left(z_l^*\big(S_l^t + jv_{\mathrm{up}(l)}b_l\big)\right) + |z_l|^2 f_l, \forall l \in \mathcal{L} \tag{8.b}$$

$$f_l = \frac{\left|S_l^t + jv_{\text{up}(l)}b_l\right|^2}{v_{\text{up}(l)}} = \frac{\left|S_l^b - jv_l b_l\right|^2}{v_l}, \quad \forall l \in \mathcal{L} \tag{8.c}$$

$$S_l^b = s_l + \sum_{m \in \mathcal{L}} \mathbf{G}_{l,m} S_m^t, \quad \forall\, l \in \mathcal{L} \tag{8.d}$$

Equations (8.a), (8.b), and (8.c) are directly derived by applying the Kirchhoff's law to Fig. 1. Equation (8.d) represents the complex power-flow of line $l$ at its bottom side (See Fig.1). $S_l^b$ is a derived variable, which is introduced here to simplify the notation.

Note that (8.c) represents the square of current flow in the central part of the two-port Π model of the line. It is worth noting that the term $f_l$ does not represent the square of the current that we can measure at the line terminals; it is indeed an internal state variable of the two-port Π model.

### *C. Original Optimal Power Flow in Radial Networks*

We can formulate an optimization problem, called OPF, with the power-flow equations shown in (8). The objective function is generally represented by a convex one, and practical examples refer to (i) nodal voltage-magnitude deviation minimization with respect to the rated value, (ii) network resistive-losses minimization, (iii) lines' current flow minimization, and (iv) cost minimization of supplied energy. Here, we consider that the objective function is the minimization of the generation cost of dispatchable units and energy imported from the transmission network (or maximization of the energy exported to the grid). It should be noted that the minimization (resp. maximization) of energy import (resp. export) from (resp. to) the grid and the total resistive-losses minimization represent an equivalent objective. Therefore, the objective function shown in (9.a) is strictly increasing in total losses or energy import from the grid. The general optimization problem is

**Original Optimal Power Flow (O-OPF)**

$$\underset{s,S,v,f}{\text{minimize}} \sum_{l \in \mathcal{L}} \Big( \mathcal{C}\big(\Re(s_l), \Im(s_l)\big) \Big) + \mathcal{C}^e(P_1^t) \tag{9.a}$$

Subject to:

$$\text{Set of Equations (8)} \tag{9.b}$$

$$v_l \le v^{\max}, \quad \forall\, l \in \mathcal{L} \tag{9.c}$$

$$v^{\min} \le v_l, \quad \forall\, l \in \mathcal{L} \tag{9.d}$$

$$\left|S_l^b\right|^2 \le I_l^{\max} v_l, \quad \forall\, l \in \mathcal{L} \tag{9.e}$$

$$\left|S_l^t\right|^2 \le I_l^{\max} v_{\text{up}(l)}, \quad \forall\, l \in \mathcal{L} \tag{9.f}$$

$$p_l^{\min} \le \Re(s_l) \le p_l^{\max}, \quad \forall\, l \in \mathcal{L} \tag{9.g}$$

$$q_l^{\min} \le \Im(s_l) \le q_l^{\max},\ \ \forall\, l \in \mathcal{L} \tag{9.h}$$

where $\mathcal{C}(.)$ is the cost function of nodal absorption (injection), $\mathcal{C}^e(P_1^t)$ is the cost function related to energy import from the grid. Both $\mathcal{C}(.)$ and $\mathcal{C}^e(.)$ are assumed to be convex, and as mentioned above, $\mathcal{C}^e(.)$ is strictly increasing. $I_l^{\max}$ and $v^{\max}/v^{\min}$ represent the square of the current limit of the lines and the maximum/minimum of square of nodal voltage-magnitudes.

In order to account for the voltage and current operational constraints, Equations (9.c)-(9.f) are added to the optimization problem. It is worth noting that the lines' ampacity limits must not be applied to $f_l$ as it does not represent the exact value of the current at its terminals. Additionally, a line-ampacity limit has to be applied to both ends of the line. The constraints (9.g) and (9.h) represent the upper and lower limits of nodal absorption. Note that the power injection, $s_l$, of a bus ($l \in \mathcal{L}$) is normally constrained to be in a pre-specified set $\mathcal{S}_l$ that is not necessarily convex. The renewable resources are normally interfaced with the grid through power electronic converters with a fixed power factor or minimum power factor requirements by the grid operators. These requirements could be modeled by adding the following linear (then convex) constraints to the optimization problem:

a) fixed power factor

$$\mathcal{S}_l = \left\{ s_l \in \mathbb{C}\ | p_l^{\min} \le\ \Re(s_l) \le p_l^{\max}, |\Im(s_l)| = \sqrt{1-\wp^2}|\Re(s_l)| \right\}$$

b) minimum power factor requirement

$$\mathcal{S}_l = \left\{ s_l \in \mathbb{C}\ | p_l^{\min} \le\ \Re(s_l) \le p_l^{\max}, |\ \Im(s_l)| \le \sqrt{1-\wp^2}|\Re(s_l)| \right\}$$

where $\wp$ represents the power factor (lead or lag).

*D. Relaxed Optimal Power Flow (R-OPF)*

O-OPF is non-convex due to Equation (8.c). However, as shown in [26], it becomes convex if we replace (8.c) by (10):

$$f_l \ge \frac{\left|S_l^t + j v_{\mathrm{up}(l)} b_l\right|^2}{v_{\mathrm{up}(l)}},\ \ \forall l \in \mathcal{L} \tag{10}$$

The new problem obtained by such a replacement is called Relaxed Optimal Power Flow (R-OPF). It can be easily shown that R-OPF is a convex problem. However, it could often occur that the optimal solution does not satisfy the original constraint (8.c), (i.e., the obtained solution has no physical meaning [29]). This could occur when the nodal upper voltage-magnitudes or lines' ampacity-limit, in case of reverse power flow[1], are binding (See [29] and part 1 of Section VII of this paper). The relaxed equation of $f$ implies that the active and reactive losses are relaxed. The relaxed losses could be interpreted as a non-negative consumption that does not exist in reality, but could be misused to relieve

[1] Note that we use the term *direct power flow* when the active and/or reactive power-flows are from bus up($l$) to bus $l$ is positive. When this term is negative, we refer to *reverse power flow*. The term applies to both real and imaginary parts independently.

the security constraints in case of large injections.

In the following, we present an augmented formulation of the O-OPF and the R-OPF that, as we prove in the following section, does not have this problem, at the expense of slight additional constraints.

## III. INTRODUCING NEW VARIABLES AND CONSTRAINTS

The main idea for modifying the OPF problem is to put the security constraints on a set of variables that (i) are upper bound for nodal voltage-magnitudes and current magnitude of the lines and (ii) do not depend on $f$. These are achieved by introducing an upper bound ($\bar{f}$) and a lower bound (a vector of zeroes) for $f$. The upper and lower bounds of $f$ are used to define the above-mentioned set of constraints. Note that the case of lower bound of $f$ (a vector of zeroes) is known in the literature as the linear DistFlow formulation. In this respect, we first introduce the following sets of auxiliary variables $\bar{f}$, $\hat{S}$ , $\bar{S}$ for the lines of the grid and $\bar{v}$ for the buses of the network, as defined in (11).

$$\hat{S}_l^t = s_l + \sum_{m\in\mathcal{L}} \mathbf{G}_{l,m}\hat{S}_m^t - j\left(\bar{v}_{\mathrm{up}(l)} + \bar{v}_l\right)b_l,\ \ \forall\, l \in \mathcal{L} \tag{11.a}$$

$$\bar{v}_l = \bar{v}_{\mathrm{up}(l)} - \ 2\Re\left(z_l^*\left(\hat{S}_l^t + j\bar{v}_{\mathrm{up}(l)}b_l\right)\right),\ \ \forall\, l \in \mathcal{L} \tag{11.b}$$

$$\bar{S}_l^t = s_l + \sum_{m\in\mathcal{L}} \mathbf{G}_{l,m}\,\bar{S}_m^t + z_l\bar{f}_l - j\left(v_{\mathrm{up}(l)} + v_l\right)b_l,\ \ \forall\, l \in \mathcal{L} \tag{11.c}$$

$$\bar{f}_l v_l \geq \max\left\{\left|\hat{P}_l^b\right|^2, \left|\bar{P}_l^b\right|^2\right\} + \max\left\{\left|\hat{Q}_l^b - \bar{v}_l b_l\right|^2, \left|\bar{Q}_l^b - v_l b_l\right|^2\right\},\ \ \forall l \in \mathcal{L} \tag{11.d}$$

$$\bar{f}_l v_{\mathrm{up}(l)} \geq \max\left\{\left|\hat{P}_l^t\right|^2, |\bar{P}_l^t|^2\right\} + \max\left\{\left|\hat{Q}_l^t + \bar{v}_{\mathrm{up}(l)}b_l\right|^2, \left|\bar{Q}_l^t + v_{\mathrm{up}(l)}b_l\right|^2\right\},\ \ \forall l \in \mathcal{L} \tag{11.e}$$

$$\bar{S}_l^b = s_l + \sum_{m\in\mathcal{L}} \mathbf{G}_{l,m}\,\bar{S}_m^t,\ \ \forall\, l \in \mathcal{L} \tag{11.f}$$

$$\hat{S}_l^b = s_l + \sum_{m\in\mathcal{L}} \mathbf{G}_{l,m}\hat{S}_m^t,\ \ \forall\, l \in \mathcal{L} \tag{11.g}$$

***Lemma I:*** *If* $(s, S, v, f, \hat{S}, \bar{v}, \bar{f}, \bar{S})$ satisfies (8) and (11), then:

1- $f \leq \bar{f}$, $v \leq \bar{v}$, $\hat{P}^t \leq P^t \leq \bar{P}^t$, and $\hat{Q}^t \leq Q^t \leq \bar{Q}^t$

2- If $(s, S, v, f)$ satisfies (8) and $\left(s, S', v', f', \hat{S}, \bar{v}, \bar{f}', \bar{S}'\right)$ satisfies (8.a), (8.b), (8.d), (10), (11) with $0 < v' \leq v$, then $\exists\ \left(\bar{f}, \bar{S}\right)$ such that $\bar{f} \leq \bar{f}'$, $\bar{P} \leq \bar{P}'$, $\bar{Q} \leq \bar{Q}'$ and $(s, S, v, f, \hat{S}, \bar{v}, \bar{f}, \bar{S})$ satisfies (8) and (11).

The proof of Lemma I is in Appendix II. Lemma I implies that $\hat{P}^t$ and $\hat{Q}^t$ represent lower bounds on $P^t$ and $Q^t$, respectively, and are adapted from linear DistFlow equations [8]. $\bar{S}$, $\bar{f}$, and $\bar{v}$ are upper bounds on $S$, $f$ and $v$, respectively.

## IV. Augmented Relaxed Optimal Power Flow

The following Augmented OPF (A-OPF) is formulated by adding the set of Equations (9) and (11), which gives the Equations (12), as follows.

**Augmented Optimal Power Flow (A-OPF)**

$$\underset{s,S,v,f,\hat{S},\hat{v},\bar{S},\bar{f}}{\text{minimize}} \sum_{l \in \mathcal{L}} \Big(\mathcal{C}\big(\Re(s_l),\Im(s_l)\big)\Big) + \mathcal{C}^e(P_1^t)) \tag{12.a}$$

subject to

$$(8), (9.g), (9.h), (11) \tag{12.b}$$

$$v^{\min} \le v_l, \quad \forall\, l \in \mathcal{L} \tag{12.c}$$

$$\bar{v}_l \le v^{max}, \quad \forall\, l \in \mathcal{L} \tag{12.d}$$

$$\Big|\big|\max\{|\hat{P}_l^b|,|\bar{P}_l^b|\}\big| + j\max\{|\hat{Q}_l^b|,|\bar{Q}_l^b|\}\Big|^2 \le v_l I_l^{\max}, \forall l \in \mathcal{L} \tag{12.e}$$

$$\big|\max\{|\hat{P}_l^t|,|\bar{P}_l^t|\} + j(\max\{|\hat{Q}_l^t|,|\bar{Q}_l^t|\}\big|^2 \le v_{\mathrm{up}(l)} I_l^{\max}, \forall l \in \mathcal{L} \tag{12.f}$$

$$P_l^t \le \bar{P}_l^t \le P_l^{\max},\; Q_l^t \le \bar{Q}_l^t \le Q_l^{\max}, \forall l \in \mathcal{L} \tag{12.g}$$

In the A-OPF, the upper limit of voltage magnitudes is imposed on $\bar{v}$, an upper bound of $v$. This constraint is shown in Equation (12.d). Similarly, the lines' current limit is modeled using the maximum of absolute values of $\bar{P}$(resp. $\bar{Q}$) and $\hat{P}$(resp. $\hat{Q}$), the upper and lower bounds of $P$ (resp. $Q$). We also add the constraint (12.g). Note that (12.g) is not a physical constraint of the system. We add it for technical ease and to more straightforwardly obtain the exactness conditions. The values of $P_l^{\max}$, and $Q_l^{\max}$ can be chosen so that these constraints do not affect the feasible solution-space of A-OPF (by performing a load flow with maximum consumption/injection level of the system and obtain the maximum possible values of $P_l/\bar{P}_l$ and $Q_l/\bar{Q}_l$). In *Lemma I* we show that $\bar{P}_l^t$ and $\bar{Q}_l^t$ are upper bounds for $P_l^t$ and $Q_l^t$, respectively. Thus, (12.g) does not shrink the feasible solution space.

***Lemma II:*** The feasible solution-space of the A-OPF is a subset of the feasible solution-space of the O-OPF.

The proof of this lemma is provided in Appendix III.

*Lemma II* states that the constraints of the A-OPF are more restrictive than the O-OPF ones. Hence, the new set of constraints (12) is more conservative, with respect to the set of equations (9), and slightly shrinks the feasible solution-space. However, the removed space covers an operation zone close to the upper bound of nodal voltage-magnitudes and lines' ampacity limits that is not a desirable operating point of the network.

The A-OPF is not convex due to Equation (8.c). We can make it convex by replacing (8.c) with (10). This gives the following proposed convex OPF problem:

**Augmented Relaxed Optimal Power Flow (AR-OPF)**

$$\underset{s,S,v,f,\hat{S},\hat{v},\bar{S},\bar{f}}{\text{minimize}} \sum_{l\in\mathcal{L}} \Big(\mathcal{C}\big(\Re(s_l),\Im(s_l)\big)\Big) + \mathcal{C}^e(P_1^t) \tag{13.a}$$

Subject to:

(8.a), (8.b), (8.d) (9.g), (9.h), (10), (11), (12.c)-(12.g)

## V. FORMULATION OF CONSTRAINTS IN MATRIX FORM

The Equations (8.a), (8.a) and (11.a) can be rewritten in matrix form as follows. Vectors such as $p, v, f, P$ and matrices such as $\mathbf{H}, \mathbf{G}, \mathbf{D}$ are defined in Section II.A.

$$\hat{P} = \mathbf{H}p \tag{14}$$

$$\hat{Q} = \mathbf{H}q - \mathbf{H}\text{diag}(b)\big(\mathbf{I} + \mathbf{G}^{\mathbf{T}}\big)\bar{v} \tag{15}$$

$$P = \hat{P} + \mathbf{H}\text{diag}(r)f \tag{16}$$

$$Q = \hat{Q} + \mathbf{H}\text{diag}(x)f + \mathbf{H}\text{diag}(b)\big(\mathbf{I} + \mathbf{G}^{\mathbf{T}}\big)\mathbf{D}f \tag{17}$$

$$v = \bar{v} - \mathbf{D}f \tag{18}$$

We are also interested in $Q_l^t + v_{\text{up}(l)}b_l$ and $\hat{Q}_l^t + \bar{v}_{\text{up}(l)}b_l$, specifically the power flows inside the longitudinal components, which we call $Q_l^c$ and $\hat{Q}_l^c$, later we use them in *Lemma III*.

$$\hat{Q}^c = \hat{Q} + \text{diag}(b)\mathbf{G}^{\mathbf{T}}\hat{v} = \mathbf{H}q - \mathbf{H}\text{diag}(B)\bar{v} \tag{19}$$

$$Q^c = Q + \text{diag}(b)\mathbf{G}^{\mathbf{T}}v = \hat{Q}^c + \mathbf{F}f \tag{20}$$

where $\mathbf{F}$ is defined in (6).
The derivation of these equations are provided in Appendix I.

## VI. EXACTNESS OF AR-OPF

In this section, we provide conditions under which the relaxation (10) in (AR-OPF) is guaranteed to be exact. They can easily be verified *ex ante* from the static parameters of the grid.

### A. Statement of the Conditions:

The five conditions are as follows (matrices $\mathbf{D}, \mathbf{E}$ and $\mathbf{H}$ are defined in (1), (2) and (7)).

*Condition C1*:

$$\|\mathbf{H}^{\mathrm{T}}\mathbf{M}\| < 1$$

*Condition C2*:

$$\|\mathbf{E}\| < 1$$

*Condition C3*: there exists an $\eta_5 < 0.5$ such that

$$\mathbf{DE} \leq \eta_5 \mathbf{D}$$

*Condition C4*: there exists an $\eta_1 < 0.5$ such that

$$(\mathbf{H}\text{diag}(r)\ \mathbf{E}) \circ \mathbf{H} \leq \eta_1 \mathbf{H}\text{diag}(r)$$

*Condition C5*: there exists an $\eta_2 < 0.5$ such that

$$(\mathbf{H}\text{diag}(r)\ \mathbf{EE}) \leq \eta_2 \mathbf{H}\text{diag}(r)\ \mathbf{E}$$

Concerning the interpretation of the above conditions, *C1* implies that $\mathbf{I} - \mathbf{G}^T - \mathbf{M}$ is invertible and has non-negative entries. *C2* ensures the convergence of the proposed iterative power-flow solution process. Condition *C3* implies that the voltage magnitude of all buses increases when one or more than one entry of $f$ decreases. Finally, *C3* and *C4* ensure that if $f$ (specifically the losses on line $l$) decreases, the direct active power-flow of all the lines upstream of line $l$ decreases.

*B. Exactness:*

***Theorem I**:*

***1)*** Under conditions C1-C3:

For every feasible solution $(s, S, v, f, \hat{S}, \bar{v}, \bar{f}, \bar{S})$ of AR-OPF there exists a feasible solution $(s, S^*, v^*, f^*, \hat{S}, \bar{v}, \bar{f}^*, \bar{S}^*)$ of A-OPF and also O-OPF with the same power-injection vector $s$.

**2)** Under conditions C1-C5:

Every optimal solution $(s, S, v, f, \hat{S}, \bar{v}, \bar{f}, \bar{S})$ of AR-OPF satisfies (8.c) and is thus an optimal solution of A-OPF.

Part 1 of *Theorem I* implies that the vector of absorptions $(s)$ of any feasible solution of the proposed OPF formulation belongs to a region where the upper and lower limits of nodal voltage-magnitudes and the lines' ampacity limits are satisfied. This is where we use *C2-C3* (C1 is related to the existence of Matrix **C**). Part 2 of *Theorem I* is the exactness of the relaxation. Here we use *C4* and *C5*.

The main idea of the proof of *Theorem I* is as follows. If $(s, S, v, f, \hat{S}, \bar{v}, \bar{f}, \bar{S})$ is feasible for AR-OPF, then $(S, v, f)$ is in general not a load-flow solution for the power injections $s$ (as (10) replaces (8.c)) but it is always possible to replace $(S, v, f)$ by $(S^*, v^*, f^*)$ obtained by a performing a load-flow on $s$. The technical difficulties are to show existence of such a load-flow solution and to find the good one (as there are multiple solutions), specifically the one that satisfies the voltage and ampacity constraints. This "good" load-flow solution is obtained using an ad-hoc iterative scheme shown in algorithm I. Furthermore, we show that an optimal solution of AR-OPF is also a load-flow solution. This is where conditions *C4* to *C5* are required. The Theorem I is proved using Lemma III introduced in the following.

*Algorithm 1:* (apexes represent iteration numbers)

Input: $\omega = \left(s, P, Q, v, f, \hat{P}, \hat{Q}, \hat{v}, \bar{f}, \bar{S}\right)$

Initialization.

$$f^{(0)} \leftarrow f$$

$$v^{(0)} \leftarrow v$$

$$P^{(0)} \leftarrow P$$

$$Q^{c(0)} \leftarrow Q^{c}$$

$$n = 1$$

for $n \geq 1$

Step 1: $f_l^{(n)} \leftarrow \frac{\left(P_l^{t(n-1)}\right)^2 + \left(Q_l^{c(n-1)}\right)^2}{v_{\mathrm{up}(l)}^{(n-1)}}$

Step 2: $P^{(n)} \leftarrow \hat{P} + \mathbf{H}\mathrm{diag}(r) f^{(n)}$ (Eq. (16))

Step 3: $Q^{c(n)} \leftarrow \hat{Q}_l^c + \mathbf{F} f^{(n)}$ (Eq. (20))

Step 4: $v^{(n)} \leftarrow \bar{v} - \mathbf{D} f^{(n)}$ (Eq. (18).)

***Lemma III:*** Under conditions *C1*-*C5*, let $\eta = \max(\eta_1, \eta_2, \eta_3, \eta_4, \eta_5) < 0.5$. For $n \geq 1$:

$$\left|\Delta f^{(n)}\right| \leq \mathbf{E}^{n-1}\left|\Delta f^{(1)}\right| \tag{21}$$

$$\left|\Delta v^{(n)}\right| \leq \eta^{n-1}\left|\Delta v^{(1)}\right| \tag{22}$$

$$v \leq v^{(n)} \tag{23}$$

$$\left|\Delta P_l^{t(n)}\right| \leq \eta^{n-1}\left|\Delta P_l^{t(1)}\right|, \forall\, l \in \mathcal{L}^{\mathcal{M}} \tag{24}$$

$$P_l^{t(n)} \leq P_l^t, \forall\, l \in \mathcal{L}^{\mathcal{M}} \tag{25}$$

$$f_l^{(n)} \leq \bar{f}_l, \forall\, l \in \mathcal{L} \tag{26}$$

$$P_l^{t(n)} \leq \bar{P}_l^t, \forall\, l \in \mathcal{L} \tag{27}$$

$$Q_l^{c(n)} \leq \bar{Q}_l^t + b_l v_{\mathrm{up}(l)}, \forall\, l \in \mathcal{L} \tag{28}$$

where $\Delta f^{(n)} = f^{(n)} - f^{(n-1)}$ for $n \geq 1$ and similarly with $P$, $Q^c$ and $v$.

The proof of Lemma III is provided in Appendix V.

## C. *Proof of Theorem I*

**Item 1:** Let $\omega = \left(s, P + jQ, v, f, \hat{P} + j\hat{Q}, \bar{v}, \bar{f}, \bar{S}\right)$ be a feasible solution of AR-OPF. Let $\mathcal{L}^{\neq}$ be the set of lines where (10) holds with strict inequality. If $\mathcal{L}^{\neq}$ is empty, $\omega$ is a load flow solution and *Theorem I* is trivially true. Assume now that $\mathcal{L}^{\neq}$ is not empty and $M$ lines have strict inequality in (10). Denote the set of the lines with strict inequality $\mathcal{M} = \{1,2,\dots,M\}$. We denote the lines for which $\mathbf{H}_{l,m} = 1, \forall\, m \in \mathcal{M}$ upstream lines/buses $\left(\mathcal{L}^{\mathcal{M}}\right)$. We now create a load flow solution for $s$. Using *Lemma III,* first we will show that, under conditions *C1-C3,* the created load flow solution is feasible (satisfies the constrains of A-OPF); then we show that, under conditions *C1-C5,* it has a lower objective function.

Consider the sequence $(P^{(n)}, Q^{c(n)}, v^{(n)}, f^{(n)})$ defined for $n \geq 0$ by means of *Algorithm I*. We now show that this sequence converges.

For $n \geq 1$ let $\Delta f^{(n)} \triangleq f^{(n)} - f^{(n-1)}$. Using *Lemma III* we have

$$\left\|\Delta f^{(n)}\right\| \leq \|\mathbf{E}\|^{(n-1)} \left\|\Delta f^{(1)}\right\| \tag{29}$$

when *C2* holds we have

$$\|\mathbf{E}\| < 1.$$

This implies that $\left\|\Delta f^{(n)}\right\| \to 0$ as $n \to \infty$, which implies that the sequence $f^{(n)}$ converges. It follows that $(P^{(n)}, Q^{c(n)}, v^{(n)}, f^{(n)})$ converges to some limit, say $(P^{*}, Q^{c*}, v^{*}, f^{*})$.

Since $f^{(n)} \geq 0$ by construction, it follows also that $f^{*} \geq 0$, and since $\mathbf{D}$ is non-negative and from item 3 of *Lemma III* (*C3* is used here. See proof of *Lemma III*).

$$v^{\min} \leq v \leq v^{*} \leq \bar{v} \leq v^{\max} \tag{30}$$

Furthermore, by step 1

$$f_l^{*} = \frac{(P_l^{t*})^2 + (Q_l^{c*})^2}{v_{\mathrm{up}(l)}^{*}} \text{ for all } l \in \mathcal{L} \qquad (31).$$

Let $Q^{*} = Q^{c*} - \mathrm{diag}(b)\mathbf{G}^T v^{*}$. It follows that $(s, S^{*} = P^{*} + jQ^{*}, v^{*}, f^{*})$ satisfies (10) with equality, i.e., it is a load flow solution and satisfies (8). Furthermore, item (2) of *Lemma I* guarantees that there exist $\bar{P}^{*}$ and $\bar{Q}^{*}$ such that $\bar{P}^{*} \leq \bar{P} \leq P^{\max}$ and $\bar{Q}^{*} \leq \bar{Q} \leq Q^{\max}$. Using $(s, S^{*}, v^{*}, \hat{S}, \bar{v}, \bar{S}^{*}, f^{*})$ and Equations (11.d) and (11.e) we can create $\bar{f}^{*}$. Let $\omega^{*} = (s, S^{*}, v^{*}, f^{*}, \hat{S}, \bar{v}, \bar{f}^{*}, \bar{S}^{*})$. From (30), we can observe that the voltage security constraints are satisfied (constraints (12.c) and (12.d). Furthermore From item (1) of *Lemma I*, (16), and (17) we have $\hat{P} \leq P^* \leq \bar{P}^* \leq P^{\max}$ and $\hat{Q} \leq Q^* \leq \bar{Q}^* \leq Q^{\max}$ which show that constraints (12.e), (12.f), and (12.g) are also satisfied. Thus $\omega^*$ is a feasible solution of AR-OPF and of A-OPF.

Furthermore, from *Lemma II*, we have that every feasible solution of A-OPF is also feasible for O-OPF.

This proves the first item of *Theorem I*.

■

**Item 2:** Assume that $\omega$ is an optimal solution of AR-OPF but not a feasible solution of A-OPF, i.e., $\mathcal{L}^{\neq}$ is non-empty and $M$ lines have strict inequality in (10). First note that at the first line ($l = 1$) we have ($\mathbf{H}_{1l} \neq 0,\ \forall\, l \in \mathcal{L}$), thus it is always in the set of upstream lines and we have

$$P_1^{t(1)} - P_1^t = \sum_{l \in \mathcal{L}^{\neq}} \mathbf{H}_{1,l}\, r_l \left(f_l^{(1)} - f_l\right) = \sum_{l \in \mathcal{L}^{\neq}} r_l \left(f_l^{(1)} - f_l\right) < 0 \tag{32}$$

, thus $P_1^{t(1)} < P_1^t$. Furthermore, by item 4 of *Lemma III* (Equation (24)) we have (*C4* and *C5* are used here. See proof of *Lemma III*):

$$\left|P_1^{t(n)} - P_1^{t(1)}\right| \leq \left|\Delta P_1^{t(n)}\right| + \cdots + \left|\Delta P_1^{t(2)}\right| \leq (\eta^{n-1} + \cdots + \eta)\left|\Delta P_1^{t(1)}\right| \leq \frac{\eta}{1-\eta}\left|\Delta P_1^{t(1)}\right| = \frac{\eta}{1-\eta}\left(P_1^t - P_l^{t(1)}\right)$$

Since $0 < \eta < 0.5$ and $P_1^{t(1)} < P_1^t$ thus

$$\left|P_1^{t(n)} - P_1^{t(1)}\right| \leq \frac{\eta}{1-\eta}\left(P_1^t - P_l^{t(1)}\right) < \left(P_1^t - P_l^{t(1)}\right) \tag{33}$$

therefore $(P_1^t)^* < P_1^t$. Now $P_1^t$ [resp. $(P_1^t)^*$ ] is the net active power import from the external grid for the solution $\omega$, [resp. $\omega^*$]. Since the power injections $s$ are identical for $\omega$ and $\omega^*$, it follows that the objective function of $\omega^*$ is strictly less than that of $\omega$, which contradicts the optimality of $\omega$.

This proves the second item of *Theorem I*.

■

To summarize, the AR-OPF is a combination of the original load-flow equations and the linear DistFlow [8] models with the inclusion of transverse parameters of the lines. Under the sufficient conditions provided above, the feasible solution-space of the AR-OPF is a subset of the one of the O-OPF, whereas the solution of R-OPF could lay outside the feasible solution-space of O-OPF. These concepts are schematically represented in Fig. 2 (see *Lemma II*).

Note how our method differs from previous relaxation-based ones. Indeed, in addition to the proper inclusion of shunt elements and lines' ampacity limits, we use supplementary variables that, as we show in *Lemma I* and (15.d)- (15.f), are bounds to the true physical quantities. Next, we require that the security constraints apply to the original, as well as to the supplementary variables. Only then we apply an SOCP relaxation of some of the constraints. A fundamental difference from the current literature shows up in the first item of Theorem 1 where we obtain that for any feasible solution of the relaxed optimization problem --not only the optimal one -- the vector of power injections corresponds to one solution of the original, non-relaxed OPF problem (i.e., it is physically feasible).

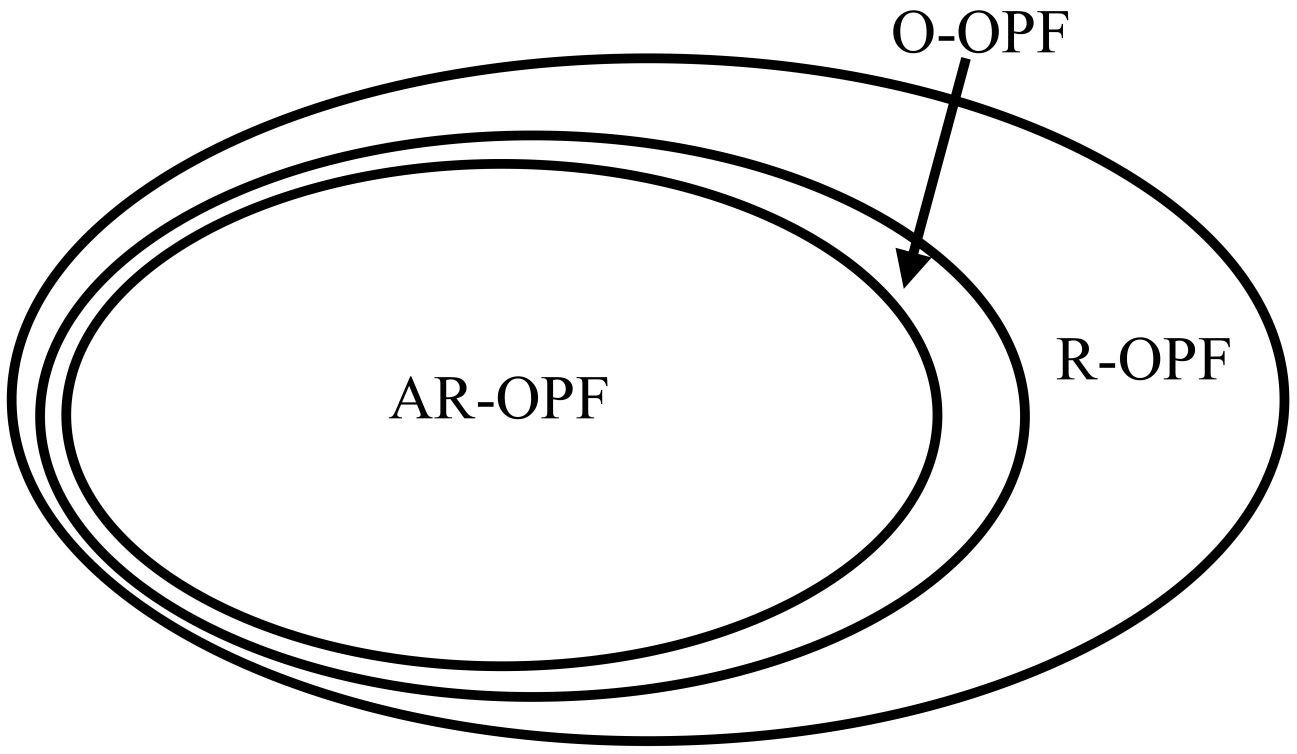

Fig. 2: Feasible solution-spaces of O-OPF, R-OPF, and AR-OPF under the five conditions provided in Section III.A

## *D. Validity of the Conditions as a Function of the Network Electrical Parameters and Physical Extension*

Note that conditions *C1-C5* are a function of the network topology and its electrical parameters. It is of interest to make observations about the validity of *C1-C5* that are functions of the grid's physical characteristics.

For a power system characterized by a given rated voltage, the per-unit-length (*pul*) electrical parameters of line, $z_l^{pul}$ and $b_l^{pul}$, do not vary drastically [31]. Also, parameters $z_l$ and $b_l$ are linear with the line lengths $\mathfrak{L}_l$.

By expressing the left-hand side of *C1* as a function of the line *pul* parameters and $\mathfrak{L}_\ell$, we note that it is given by $\left(\max_{l\in\mathcal{L}} x_l^{pul}\right)\left(\max_{l\in\mathcal{L}} B_l^{pul}\right)\mathfrak{L}_l^2$. Also for *C2*, it is straightforward to observe that the $l_1$–norm of matrix $\mathbf{E}$ has a linear dependency with $\mathfrak{L}_l$. Concerning *C3*, *C4*, and *C5* we can observe that the left-hand side of their inequality is proportional to $\mathfrak{L}_l^2$, whereas the right-hand side of their inequality is proportional to $\mathfrak{L}_l$.

Therefore, for given *pul* parameters $z_l^{pul}$ and $b_l^{pul}$, there exists a $\mathfrak{L}_l$ small enough that *C1-C5* holds. The consequence of this observation is that *C1-C5* can be verified a priori for families of networks characterized by given electrical parameters and physical extensions. In the following section, we numerically show that conditions *C1-C5* hold, with large margins, for the tested real distribution networks.

# VII. Simulation results

This section is divided into three parts. The first part demonstrates the capabilities of the proposed model to provide feasible solutions as well as infeasible behaviors of the existing convex OPF models. The influence of the inclusion of the shunt elements is also discussed. In the second part, the IEEE 34-bus network [32] and CIGRE European benchmark medium-voltage network [33] are selected to assess the scalability of the provided conditions. Finally, in the third part it is shown that the compression of solution space associated with the conservative constraints of the AR-OPF is small.

### A. Comparison of AR-OPF with R-OPF and AR-OPF without Shunt Elements

The simple network introduced in [29] is chosen to show the capability of the AR-OPF to provide an optimal and feasible solution. The grid is composed of three identical coaxial power cables. Fig. 3 shows the topology of the grid. The cable data is presented in Table I. Note that the values of the line parameters refer to the typical underground cables in use in actual distribution networks[1].

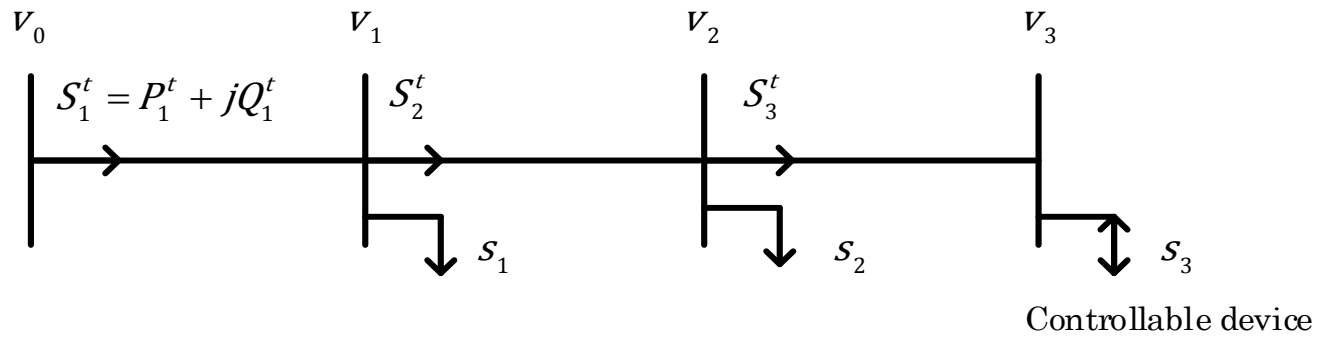

Fig. 3: Network used for comparison of different models

TABLE I. NETWORK PARAMETERS.

| Parameter | Value |
|---|---|
| Line parameters $R\left(\frac{\text{ohm}}{\text{km}}\right), L\left(\frac{\text{mH}}{\text{km}}\right), C\left(\frac{\mu\text{F}}{\text{km}}\right)$, length(km) | (0.193, 0.38, 0.24, 1) |
| Network rated voltage and base power (kV, MVA) | (24.9, 5) |
| Power rating (MW) (Storage on bus 3) | 1.5 |
| $(v^{min}, v^{max})$ (p.u.) | $(0.9 \times 0.9, 1.1 \times 1.1)$ |
| $I_{kl}^m$ (for all 3 lines) ($A^2$) | 120× 120 |
| Complex load (3 phase) on bus 2 and 3 (p.u.) | $(-0.21 - j0.126)$, $(-0.252 - j0.1134)$ |
| Energy cost from external grid, cost of active power production/consumption of controllable device at bus 3 ($/MWh) | (150, 50) |

The objective is to minimize the total cost of imported electricity, plus the cost of active power production/consumption of the controllable device connected to bus 3, assumed to be an energy storage system. Table I contains all the elements considered in the cost function. The formulated AR-OPF is the following:

$$\underset{s,S,v,f,\hat{S},\hat{v},\bar{S},\bar{f}}{\text{minimize}} \; 50\Re(s_3) + 150P_1^t$$

Subject to:

[1] They are derived from page 16 of the following document. http://www.nexans.com/Switzerland/files/NEXANS06_BTMTAcc_F.pdf

(8.a), (8.b), (8.d) (9.g), (9.h), (10), (11), (12.c)-(12.h)

$$s_1 = -0.21 - j0.126$$
$$s_2 = -0.252 - j0.1134$$
$$-0.3 \le \Re(s_3) \le 0.3$$
$$\Im(s_3) = 0$$

The lines' current magnitudes for three cases AR-OPF, R-OPF, and a case where shunt elements are not considered in AR-OPF are shown in Fig. 4. In particular, Fig. 4-a shows the current flow of the lines for the solution obtained with the AR-OPF. Fig. 4-b and 4-c correspond to R-OPF and the case where shunt elements are not considered in AR-OPF, respectively. These current flows are calculated using an a posteriori load-flow analysis. It can be seen that the maximum rating of the lines (dashed line in Fig. 4) is satisfied only with the solution provided by the AR-OPF, whereas they are largely violated in the two other cases.

*B. Scalability of the Conditions for Benchmark Networks*

The scalability analysis is done by increasing the maximum level of injections into the systems. We choose these two grids because the former network is composed by long overhead lines, whereas underground cables with high penetration of distributed energy resources characterize the latter. Both networks are considered to be balanced. The minimum and maximum nodal voltage-magnitude limits are considered to be 0.95 and 1.05 p.u., respectively.

For the IEEE 34 buses, the line impedances are the positive sequences (it is assumed that the gird is a three-phase balanced one). The base apparent-power and voltage-magnitude values are chosen to be 5 MVA and 24.9 kV, respectively. We increase the active power-injections at each bus proportionally to their load share, because there are no DGs in the network except for the two shunt capacitors. $P_l^{\max}$ and $Q_l^{\max}$ for each line is considered to be 110% of total active and reactive load in their downstream nodes, respectively. The first condition that is violated is C3. However, this condition is violated for a total net injection equal to 2.2 MW. For this operating point, the nodal voltage-magnitudes reach a maximum value of 1.073 p.u., a value far from typical feasible operating conditions.

For the CIGRE network, the positive-sequence impedance, with base apparent-power and voltage-base values equal to 25 MVA and 20 kV are used. The network already has a 3.079 MW generation capacity (it is designed to study the penetration of renewable resources). $P_l^{\max}$ and $Q_l^{\max}$ for each line is considered to be 110% of total active and reactive load in their downstream nodes, respectively. The first condition that is violated is C3. These violation occur at 585% of DG penetration, corresponding to 18.01 MW of active power production. For this operating point, the maximum value of the nodal voltages reaches 1.0857 p.u., again a value far from feasible.

Concerning the applicability of the proposed formulation to unbalanced systems, in the case of electrically balanced systems, we can apply the proposed AR-OPF to each sequence separately. For electrically unbalanced systems, work is in progress to extend the proposed approach.

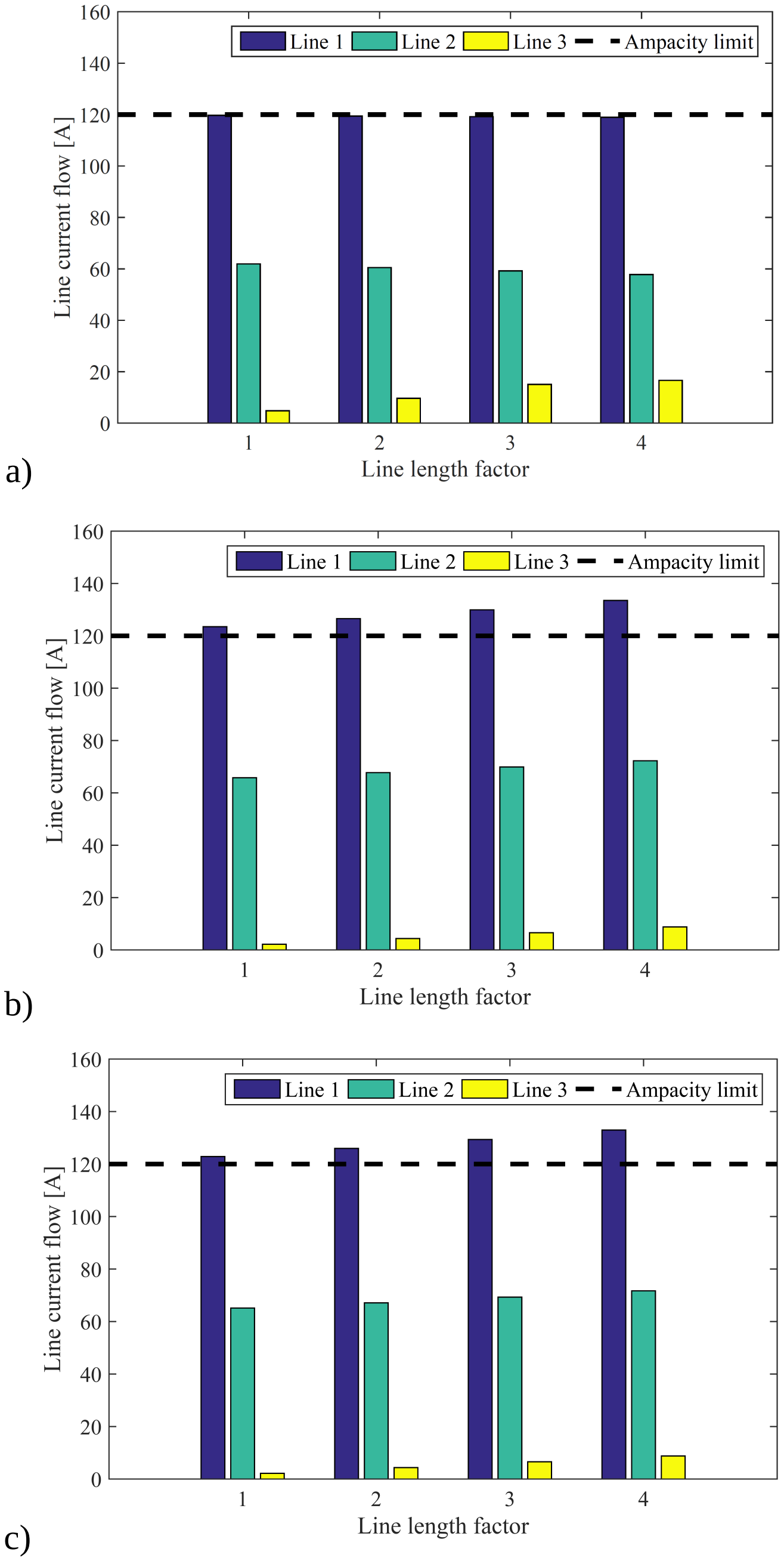

Fig. 4: Current flow magnitude of the lines vs. cables length. (a) AR-OPF, (b) R-OPF, and (c) AR-OPF without transverse parameters (nominal voltage 24.9 kV))

### *C. Quantification of the Compression of the Solution Space Associated with the Conservative Constraints of AR-OPF*

For the two case studies reported in this paper, we numerically show that the compression of the solution space is small (for this analysis we used a scenario with high level of nodal power injections). These analyses are reported in the following.

#### *1) IEEE 34-bus Network:*

*Nodal voltage-magnitudes:*

The step-by-step procedure to carry out this analysis is reported here.

First, we relax the ampacity limits of the lines (Equations (12.e) and (12.f)). Then, we increase the nodal injections up to the point where the constraint (12.d) ($\hat{v}_l \leq v^{max} = (1.05)^2$) is binding for at least one of the buses (note that the value of 1.05 p.u. is a typical upper bound for the nodal voltage-magnitudes imposed by power quality norms).

By analyzing the results we find that, for this network, the first binding voltage-magnitude corresponds to node #20. The maximum difference between the nodal voltage-magnitudes ($|V_l| = \sqrt[2]{v_l}$) and the corresponding auxiliary ones $\left(|\bar{V}_l| = \sqrt[2]{\bar{v}_l}\right)$ is equal to 0.0011 ($|\bar{V}_l| = 1.05$ p.u., $|V_l| = 1.048943$ p.u.). This value supports the claim about the small compression of the solution space for the node voltage-magnitudes.

*Lines' ampacity limits:*

The same analysis is done here for lines' ampacity limits. This time the nodal voltage-magnitudes are relaxed (Equations (12.c) and (12.d)). Note that the lines' current flows are bounded by Equations (9.e) and (9.f), in the O-OPF. The corresponding constraints in the AR-OPF are (12.e) and (12.f). We increase the nodal injections to the point where at least one of the lines' ampacity limits becomes binding ((12.e) and/or (12.f)). In this analysis, line # 1 is the first one (its current limit is 1 p.u. where base current value is 104.34 A). For the obtained operating point, the auxiliary current magnitude $\left(\frac{\left|\max\{|\hat{P}_l^b|,|\bar{P}_l^b|\}+j\max\{|\hat{Q}_l^b|,|\bar{Q}_l^b|\}\right|^2}{v_l}\right)$ is equal to 1 p.u., whereas the original current flow $\left(\frac{|S_l^b|^2}{v_l}\right)$ is 0.91 p.u. (9% difference). This value supports the claim about the small compression of the solution space for the line currents constraints.

*2) CIGRE European Network:*

The same analysis, as the one for IEEE 34-bus network, is carried out for the CIGRE network.

*Nodal voltage-magnitudes:*

We have repeated the same process described for the IEEE benchmark. The first binding voltage-magnitude corresponds to the node #8. The maximum difference between the nodal voltage-magnitudes ($|V_l| = \sqrt[2]{v_l}$) and the corresponding auxiliary ones ($|\bar{V}_l| = \sqrt[2]{\bar{v}_l}$) is equal to 0.0037 p.u.. Also for this benchmark grid, this obtained value supports the claim about the small compression of the solution space for the nodal voltage constraints.

*Lines' ampacity limits:*

Similarly to the previous cases, we increased the nodal injections to the point where at least one of the lines' ampacity limits becomes binding (12.e) and/or (12.f). Here again, line # 1 is the first line that becomes binding (its current limit is 1 p.u. where base current value is 974.2786 A). At this operating point, the auxiliary current magnitude $\left(\frac{\left|\max\{|\hat{P}_l^b|,|\bar{P}_l^b|\}+j\max\{|\hat{Q}_l^b|,|\bar{Q}_l^b|\}\right|^2}{v_l}\right)$ is equal to 1 p.u., whereas the original current flow $\left(\frac{|S_l^b|^2}{v_l}\right)$ is 0.912 p.u. (9.88% difference). Also for this benchmark grid, this obtained value supports the claim about the small compression of the solution space for the line current constraints.

## VIII. DISCUSSION ON THE EXTENSION OF THE MODEL TO UNBALANCED RADIAL GRIDS

In this paper, we targeted radial power-grids operating in balanced conditions. Future works include the extension to unbalanced three-phase systems. As mentioned in the Section I, several recent papers propose the use of convex relaxation to define OPF problems in grid operating in generic unbalanced conditions [20]-[23]. However, the following main challenges still remain to be addressed in unbalanced systems: (i) proper inclusion of static security constraints (i.e., voltage-magnitude's limits and lines' current flow) and (ii) exactness of the adopted relaxations. A potential way to address these challenges by means of the proposed AR-OPF is described in the following.

Let us consider a radial power grid whose generic components connected between two of its buses are characterized by circulant shunt admittance and longitudinal impedance matrices (i.e., for a matrix of rank $n$, its eigenvectors are composed by the roots of unity of order $n$). For these grids, it is possible to decompose all the nodal/flow voltages, currents and powers with the well-known sequence transformation. The result of this transformation is composed of three symmetrical and balanced three-phase circuits, for which the SOCP relaxation we have proposed can be applied as is. The main problem with this approach, however, is the transformation of the voltage/current constraints from the phase domain to the corresponding ones in the sequence domain. Indeed, such a transformation couples the voltage/currents constraints in the sequence domain. However, it is possible to separately bind the zero and negative sequence-terms of nodal voltage-magnitude and lines' current flows by using more conservative constraints as the magnitude of these quantities are restricted by standards/norms (i.e., their maximum magnitudes are known a-priori). The binding of the zero and negative sequences associated with the voltages and currents should enable the positive sequence to be decoupled. Then, we can apply the proposed SOCP relaxation to the three sequences for which we may derive different voltage/current inequalities. Once the three problems are solved, we can transform the obtained voltage/currents/powers in the sequence domain back to the (unbalanced) phase domain.

## IX. CONCLUSIONS AND PERSPECTIVES

The OPF problem in radial distribution networks is a timely research topic driven by the need to provide a robust mathematical tool to several problems associated with the vast connection of DGs in ADNs. To solve this problem, the recent literature has discussed the adoption of the SOCP relaxation. However, this approach might provide technically infeasible solutions, depending on the flows of the powers in the lines and the inclusion of the line transverse parameters.

In order to overcome these limitations, we have formulated the AC-OPF by using the line two-port Π model. Additionally, we have augmented the problem constraints to incorporate the lines' ampacity limits. In order to preserve the exactness of the relaxed problem, we have added a new set of conservative constraints for the ampacity limits of the lines, as well as the nodal voltage-magnitudes upper bound to the problem. This set of new constraints shrinks the feasible space of the solution. Furthermore, we have provided sufficient conditions for the feasibility and optimality of the proposed OPF formulation. Using IEEE and CIGRE radial-grid benchmarks, we have shown that these conditions are mild and hold for practical distribution networks in the feasible operation ranges. In order to analyze the performances of the proposed formulation, we have used a simple example. Underground coaxial power-cables, whose transverse parameters cannot be neglected, are the lines in this simple and replicable network. Consequently, the modeling of the line ampacity constraints needs to properly account for the line transverse parameters. Using this simple network, we have showed the capabilities of the proposed model to provide feasible solutions as well as infeasible behaviors of the existing convex OPF models.

Finally, we have provided the proof of the exactness of the AR-OPF, as well as the derivation of the sufficient conditions.

In this paper, we have targeted balanced radial grids. Future works will extend the procedure presented here to the case of unbalanced multi-phase radial grids.

# X. Appendices

## *A. Appendix I: Formulation of Constraints in Matrix Form*

In this Appendix, we derive the matrix-from of the power flow equations.

For $l \neq 1$ the upstream bus of $l$, $\mathrm{up}(l)$, is the unique $j \in \{1, \dots, L\}$ such that $\mathbf{G}_{j,l} = 1$, and the voltage $v_{\mathrm{up}(l)}$ at the upstream bus of $l$ is given by $v_{\mathrm{up}(l)} = \sum_j \mathbf{G}_{j,l} v_j$, namely $v_{\mathrm{up}(l)} = (\mathbf{G}^T v)_l$.

Using Equation (8.b), we can rewrite the nodal voltage equation for $l = 1, \dots, L$ as follows (recall that $S = P + jQ$ represents the vector of $S_l^t = P_l^t + jQ_l^t$ for all lines):

$$v = \mathbf{G}^T v + v_0 e - 2\mathrm{diag}(r)P - 2\mathrm{diag}(x)Q - 2\mathrm{diag}(x)\mathrm{diag}(b)\mathbf{G}^T v + \mathrm{diag}(|z|^2)\, f \quad (34)$$

where $e = (1,0, \dots 0)^T$. Similarly (8.a) gives

$$P = p + \mathbf{G}P + \mathrm{diag}(r)f \quad (35.a)$$

$$Q = q + \mathbf{G}Q + \mathrm{diag}(x)f - \mathrm{diag}(b)\left(\mathbf{I} + \mathbf{G}^{\mathbf{T}}\right)v \quad (35.b).$$

Using (1) we can rewrite (35.a) and (35.b) as

$$P = \mathbf{H}p + \mathbf{H}\mathrm{diag}(r)f \quad (36.a)$$

$$Q = \mathbf{H}q + \mathbf{H}\mathrm{diag}(x)f - \mathbf{H}\mathrm{diag}(b)\left(\mathbf{I} + \mathbf{G}^{\mathbf{T}}\right)v \quad (36.b).$$

Similarly from (11.c) we have

$$\bar{P} = \mathbf{H}p + \mathbf{H}\mathrm{diag}(r)\bar{f} \quad (37.a)$$

$$\bar{Q} = \mathbf{H}q + \mathbf{H}\mathrm{diag}(x)\bar{f} - \mathbf{H}\mathrm{diag}(b)\left(\mathbf{I} + \mathbf{G}^{\mathbf{T}}\right)v \quad (37.b).$$

We can eliminate $P, Q$ from (34), using (36.a), (36.b) and obtain

$$\begin{aligned}\left(\mathbf{I}-\mathbf{G}^T+2\mathrm{diag}(x)\mathrm{diag}(b)\mathbf{G}^T-2\mathrm{diag}(x)\mathbf{H}\mathrm{diag}(b)(\mathbf{I}+\mathbf{G}^T)\right)v \\ = v_0 e - 2\mathrm{diag}(r)\mathbf{H}p - 2\mathrm{diag}(x)\mathbf{H}q \\ +\left(-2\mathrm{diag}(r)\mathbf{H}\mathrm{diag}(r)-2\mathrm{diag}(x)\mathbf{H}\mathrm{diag}(x)+\mathrm{diag}(|z|^2)\right)f\end{aligned} \quad (38).$$

Now use (1) and the following equation

$$\mathbf{G}\,\mathrm{diag}(b)\,\mathbf{G}^T = \mathrm{diag}(B)-\mathrm{diag}(b) \quad (39)$$

which gives

$$\begin{aligned}[\mathbf{I}-\mathbf{G}^\mathbf{T}-\mathbf{M}]v = v_0 e - 2\mathrm{diag}(r)(\mathbf{H}p+(\mathbf{H}-\mathbf{I})\mathrm{diag}(r)f) - 2\mathrm{diag}(x)(\mathbf{H}q+(\mathbf{H}-\mathbf{I})\mathrm{diag}(x)f) \\ - \mathrm{diag}(|z|^2)f\end{aligned} \quad (40)$$

where $\mathbf{M}=2\mathrm{diag}(x)\mathbf{H}\mathrm{diag}(B)$. Under condition *C1* (given in Section VI), we prove in Appendix IV that $\mathbf{C}=[\mathbf{I}-\mathbf{G}^\mathbf{T}-\mathbf{M}]^{-1}$ exists and has non-negative entries. It follows that we can solve for $v$ in Equation (40) as follows:

$$v = v_0\mathbf{C}e - 2\mathbf{C}\mathrm{diag}(r)(\mathbf{H}p) - 2\mathbf{C}\mathrm{diag}(x)(\mathbf{H}q) - \mathbf{D}f \quad (41)$$

where $\mathbf{D}$ is defined in (2). Note that $\mathbf{H}-\mathbf{I} \geq 0$, therefore $\mathbf{D}$ is non-negative. Next, we can write (11.b) as follows:

$$\bar{v} = v_0\mathbf{C}e - 2\mathbf{C}\mathrm{diag}(r)(\mathbf{H}p) - 2\mathbf{C}\mathrm{diag}(x)(\mathbf{H}q) \quad (42)$$

thus

$$v = \bar{v} - \mathbf{D}f \quad (43).$$

Since $\mathbf{D}$ and $f$ are non-negative, it follows that:

$$v \leq \bar{v} \quad (44).$$

■

The Equations (8.a) and (11.a) can be rewritten in matrix form as follows:

$$\hat{P} = \mathbf{H}p \quad (45)$$

$$\hat{Q} = \mathbf{H}q - \mathbf{H}\text{diag}(b)\big(\mathbf{I} + \mathbf{G}^{\mathbf{T}}\big)\bar{v} \tag{46}$$

$$P = \hat{P} + \mathbf{H}\text{diag}(r)f \tag{47}$$

$$Q = \hat{Q} + \mathbf{H}\text{diag}(x)f + \mathbf{H}\text{diag}(b)\big(\mathbf{I} + \mathbf{G}^{\mathbf{T}}\big)\mathbf{D}f \tag{48}$$

Since $r, b, x, \mathbf{D}, \mathbf{G}, \mathbf{H}$ and $f$ are all non-negative, thus

$$\hat{P} \leq P \tag{49}$$

$$\hat{Q} \leq Q \tag{50}$$

We are also interested in $Q_l^t + v_{\text{up}(l)}b_l$ and $\hat{Q}_l^t + \bar{v}_{\text{up}(l)}b_l$, namely the power flows inside the longitudinal components, which we call them $Q_l^c$ and $\hat{Q}_l^c$, later we use them in *Lemma III*. Using (46), (48) and (36.b) we have

$$\hat{Q}^c = \hat{Q} + \text{diag}(b)\mathbf{G}^{\mathbf{T}}\hat{v} = \mathbf{H}q - \mathbf{H}\text{diag}(B)\bar{v} \tag{51}$$

$$Q^c = Q + \text{diag}(b)\mathbf{G}^{\mathbf{T}}v = \mathbf{H}q + \mathbf{H}\text{diag}(x)f - \mathbf{H}\text{diag}(B)v = \hat{Q}^c + \mathbf{F}f \tag{52}$$

where **F** is defined in (6). Similar to (50) we have

$$\hat{Q} \leq \hat{Q}^c \leq Q^c \tag{53}$$

*B. Appendix II: Proof of Lemma I*

We prove this lemma by induction starting from the leaves of the grid. Formally, for a bus $l$, let $\text{height}(l)$ denotes its height in the tree, defined by $\text{height}(l) = 0$ when $l$ is a leaf and $\text{height}(l) = 1 + \max_{k:\, \text{up}(k)=l} \text{height}(k)$ otherwise.

**I. Base case (height =0)**

For the base case we show that *Lemma I* holds for the leaves of the system.

**a)** Suppose bus $l$ is a leaf of the network ($l \in \underline{\mathcal{L}}$) with $s_l = p_l + jq_l$ as its total load (See Fig. 1). Since $(s, S, v, f, \hat{S}, \bar{v}, \bar{f}, \bar{S})$ satisfies (8), (11), from (8.c) we have

$$f_l = \frac{(p_l)^2 + (q_l - v_l b_l)^2}{v_l}, \quad \forall\, l \in \underline{\mathcal{L}} \tag{54}$$

Since $0 < v_l \leq \hat{v}_l$, thus from (11.d)

$$0 \leq f_l \leq \frac{\max\{|(p_l)|^2, |p_l|^2\}}{v_l} + \frac{\max\{|q_l - v_l b_l|^2, |q_l - \bar{v}_l b_l|^2\}}{v_l} \leq \bar{f}_l, \qquad \forall\, l \in \underline{\mathcal{L}} \tag{55}$$

combining with (36), (37), (47), and (48) it comes that

$$\hat{P}_l^t \leq P_l^t \leq \bar{P}_l^t, \qquad \forall\, l \in \underline{\mathcal{L}}$$

$$\hat{Q}_l^t \leq Q_l^t \leq \bar{Q}_l^t, \qquad \forall\, l \in \underline{\mathcal{L}}$$

this shows item 1 of *Lemma I*.

**b)** One can choose $\bar{f}_l$ as followings so that $\bar{f}_l \leq \bar{f}_l'$ and $\bar{f}_l$ satisfy (11.d) and (11.e) (recall that $0 < v' \leq v \leq \bar{v}$):

$$\bar{f}_l' \geq \bar{f}_l = \max\left\{\begin{array}{c}\left(\dfrac{\max\{|(p_l)|^2, |p_l|^2\}}{v_l'} + \dfrac{\max\{|q_l - v_l' b_l|^2, |q_l - \bar{v}_l b_l|^2\}}{v_l'}\right), \\ \left(\begin{array}{c}\dfrac{\max\left\{\left|(p_l + r_l \bar{f}_l')\right|^2, |p_l|^2\right\}}{v_{\mathrm{up}(l)}'} + \\ \dfrac{\max\left\{\left|q_l + x_l \bar{f}_l' - v_l' b_l\right|^2, |q_l - \bar{v}_l b_l|^2\right\}}{v_{\mathrm{up}(l)}'}\end{array}\right)\end{array}\right\}$$

$$\geq \max\left\{\begin{array}{c}\left(\begin{array}{c}\dfrac{\max\{|(p_l)|^2, |p_l|^2\}}{v_l} + \\ \dfrac{\max\{|q_l - v_l b_l|^2, |q_l - \bar{v}_l b_l|^2\}}{v_l}\end{array}\right), \\ \left(\begin{array}{c}\dfrac{\max\left\{\left|(p_l + r_l \bar{f}_l'')\right|^2, |p_l|^2\right\}}{v_{\mathrm{up}(l)}} + \\ \dfrac{\max\left\{\left|q_l + x_l \bar{f}_l'' - v_l b_l\right|^2, |q_l - \bar{v}_l b_l|^2\right\}}{v_{\mathrm{up}(l)}}\end{array}\right)\end{array}\right\}, \forall\, l \in \underline{\mathcal{L}}$$

hence $\bar{f}_l \leq \bar{f}_l'$ and satisfies (11.d) and (11.e). Consequently from (37) we have $(\bar{P}_l^t) \leq (\bar{P}_l^t)'$ and $(\bar{Q}_l^t) \leq (\bar{Q}_l^t)'$. These show item 2 of *Lemma I*.

**2- Induction step**

Assume the statements in *Lemma I* are true for all buses of height $\leq n$, for some $n \geq 0$. We now show that it holds for all buses of height $\leq n + 1$. Let $k$ be a bus with $\mathrm{height}(k) = n + 1$.

**a)** For all downstream buses $l$ of $k$, we have $\mathrm{height}(l) \leq n$, thus $S_l^t \leq \bar{S}_l^t$. Furthermore from Equations (8.d), (11.f), and (11.g) it comes

$$\hat{P}_k^b \leq P_k^b \leq \bar{P}_k^b \tag{56.a}$$

$$\hat{Q}_k^b \leq Q_k^b \leq \bar{Q}_k^b \tag{56.b}$$

thus (recall that $0 < v \leq \hat{v}$)

$$f_k = \frac{|P_k^b|^2 + |Q_k^b - v_k b_k|^2}{v_k} \le \frac{\max\{|\hat{P}_k^b|^2, |\bar{P}_k^b|^2\}}{v_k} + \frac{\max\{|\hat{Q}_k^b - \bar{v}_k b_k|^2, |\bar{Q}_k^b - v_k b_k|^2\}}{v_k} \le \bar{f}_k$$

combining with (36), (37), (47), and (48) it comes that

$$\hat{P}_k^t \le P_k^t \le \bar{P}_k^t$$
$$\hat{Q}_k^t \le Q_k^t \le \bar{Q}_k^t.$$

This show item 1 of *Lemma I*.

**b)** Based on the induction assumption and Equation (11.f), we have $\hat{P}_k^b \le (\bar{P}_k^b) \le (\bar{P}_k^b)'$ and $\hat{Q}_k^b \le (\bar{Q}_k^b) \le (\bar{Q}_k^b)'$ (recall that $0 < v' \le v \le \bar{v}$). Therefore one can choose $\bar{f}_k$ as follows so that $\bar{f}_k \le \bar{f}_k'$ and $\bar{f}_k$ satisfies (11.d) and (11.e):

$$\bar{f}_k' \ge \bar{f}_k = \max\left\{ \begin{pmatrix} \dfrac{\max\{|\hat{P}_k^b|^2, |(\bar{P}_k^b)'|^2\}}{v_k'} + \\ \dfrac{\max\{|\hat{Q}_k^b - \bar{v}_k b_k|^2, |(\bar{Q}_k^b)' - v_k' b_k|^2\}}{v_k'} \end{pmatrix}, \begin{pmatrix} \dfrac{\max\{|\hat{P}_k^b|^2, |(\bar{P}_k^b)' + r_k \bar{f}_k'|^2\}}{v_{\mathrm{up}(k)}'} + \\ \dfrac{\max\{|\hat{Q}_k^b - \bar{v}_k b_k|^2, |(\bar{Q}_k^b)' + x_k \bar{f}_k' - v_k' b_k|^2\}}{v_{\mathrm{up}(k)}'} \end{pmatrix} \right\}$$

$$\ge \max\left\{ \begin{pmatrix} \dfrac{\max\{|\hat{P}_k^b|^2, |(\bar{P}_k^b)|^2\}}{v_k} + \\ \dfrac{\max\{|\hat{Q}_k^b - \bar{v}_k b_k|^2, |(\bar{Q}_k^b) - v_k b_k|^2\}}{v_k} \end{pmatrix}, \begin{pmatrix} \dfrac{\max\{|\hat{P}_k^b|^2, |(\bar{P}_k^b) + r_k \bar{f}_k|^2\}}{v_{\mathrm{up}(k)}} + \\ \dfrac{\max\{|\hat{Q}_k^b - \bar{v}_k b_k|^2, |(\bar{Q}_k^b) + x_k \bar{f}_k - v_k b_k|^2\}}{v_{\mathrm{up}(k)}} \end{pmatrix} \right\}.$$

Thus $\bar{f}_k \le \bar{f}_k'$ and $\bar{f}_k$ satisfy (11.d) and (11.e). Consequently from (37) we have $(\bar{P}_k^t) \le (\bar{P}_k^t)'$ and $(\bar{Q}_k^t) \le (\bar{Q}_k^t)'$. These show item 2 of *Lemma I*.

Both basis and inductive steps are proved, which completes the proof of *Lemma I*.

■

*C. Appendix III: Proof of Lemma II*

A-OPF contains all the constraints of O-OPF except (9.c), (9.e), and (9.f). These constraints are replaced by (12.d), (12.e), and (12.f). It suffices to show that (12.d), (12.e), and (12.f) are more restrictive than (9.c), (9.e), and (9.f). The right-hand side of the constraints are the same ((9.c) with (12.d), (9.e) with (12.e), (9.f) with (12.f)). We just need to show that the left-hand sides of relevant constraints in (12) are upper bound for those in (9).

From Lemma I, (49), and (50), we have

$$\hat{P}_k^t \leq P_k^t \leq \bar{P}_k^t$$
$$\hat{Q}_k^t \leq Q_k^t \leq \bar{Q}_k^t$$

combined with Equation (8.d) it comes that

$$\hat{P}_k^b \leq P_k^b \leq \bar{P}_k^b$$
$$\hat{Q}_k^b \leq Q_k^b \leq \bar{Q}_k^b$$

thus

$$\left||\max\{|\hat{P}_l^b|, |\bar{P}_l^b|\}| + j\max\{|\hat{Q}_l^b|, |\bar{Q}_l^b|\}\right|^2 \geq |S_l^b|^2, \forall\, l \in \mathcal{L}$$
$$|\max\{|\hat{P}_l^t|, |\bar{P}_l^t|\} + j(\max\{|\hat{Q}_l^t|, |\bar{Q}_l^t|\}|^2 \geq |S_l^t|^2, , \forall\, l \in \mathcal{L}.$$

Furthermore from (44) we have $v \leq \bar{v}$.

■

*D. Appendix IV*

In this Appendix, we prove that when *C1* holds, $\mathbf{I} - \mathbf{G}^T - \mathbf{M}$ is invertible and has non-negative entries.

We can rewrite $\mathbf{I} - \mathbf{G}^T - \mathbf{M}$ as follows:

$$\mathbf{I} - \mathbf{G}^T - \mathbf{M} = (\mathbf{I} - \mathbf{G}^T)[\mathbf{I} - (\mathbf{I} - \mathbf{G}^T)^{-1}\mathbf{M}] = (\mathbf{I} - \mathbf{G}^T)[\mathbf{I} - \mathbf{H}^{\mathrm{T}}\mathbf{M}] \quad (57).$$

We now use the identity

$$\left(\mathbf{I} - \mathbf{H}^{\mathrm{T}}\mathbf{M}\right)\left[\mathbf{I} + \mathbf{H}^{\mathrm{T}}\mathbf{M} + \left(\mathbf{H}^{\mathrm{T}}\mathbf{M}\right)^2 + \cdots\right] = \mathbf{I} \quad (58)$$

which holds whenever $||\mathbf{H}^T\mathbf{M}|| < 1$ (recall that $||\mathbf{H}^T\mathbf{M}||$ is the Frobenius norm).

It follows that, when *C1* holds, then $||\mathbf{H}^T\mathbf{M}|| < 1$, which by Equation (58) proves that $\left(\mathbf{I} - \mathbf{H}^{\mathrm{T}}\mathbf{M}\right)$ is invertible. By transposition of Equation (1), $\left(\mathbf{I} - \mathbf{G}^{\mathrm{T}}\right)$ is also invertible; together with Equation (57), this shows that $\mathbf{I} - \mathbf{G}^T - \mathbf{M}$ is invertible when *C1* holds. Furthermore, (1), (57) and (58) imply that

$$(\mathbf{I}-\mathbf{G}^T-\mathbf{M})^{-1}=\left(\mathbf{I}+\mathbf{H}^{\mathrm{T}}\mathbf{M}+\left(\mathbf{H}^{\mathrm{T}}\mathbf{M}\right)^2+\cdots\right)\mathbf{H}^T\geq 0 \quad (59).$$

■

*E. Appendix V: Proof of Lemma III*

The proof is by induction on $n\geq 1$.

**1- Base case ($n=1$):**

(21), (22), and (24) are trivially true.

We have $f_l^{(1)}\leq f_l^{(0)}$ for every $l$ because $f^{(1)}$ is the right-hand side of (10) in the original formulation of the constraints and $\omega$ is feasible. By Equations (47) and (43) since $\mathbf{D}\geq\mathbf{0}$ and $\mathbf{H}\mathrm{diag}(r)\geq 0$, it follows that $P^{(1)}\leq P$ and $v\leq v^{(1)}$. This shows (23) and (25).

Since $P^{(0)}$ and $Q^{(0)}$ are feasible solution of AR-OPF, we have:

$$\hat{P}_l^t\leq P_l^{t(0)}\leq\bar{P}_l^t\,,\forall\, l\in\mathcal{L} \quad (60)$$

$$\hat{Q}_l^c\leq Q_l^{c(0)}=Q_l^{t(0)}+b_l v_{\mathrm{up}(l)}^{(0)}\leq\bar{Q}_l^t+b_l v_{\mathrm{up}(l)}\,,\quad\forall\, l\in\mathcal{L}$$

Thus from step 1 of *Algorithm I*, Equation (11.e) and noting that $v\leq v^{(1)}$, we have $f_l^{(1)}\leq\bar{f}_l\;\forall\, l\in\mathcal{L}$. This shows Equation (26). Furthermore, knowing that $v\leq v^{(1)}$ and using Equations (36), (37) and (52) one can show that

$$P_l^{t(1)}\leq\bar{P}_l^t,\forall\, l\in\mathcal{L}$$

$$Q_l^{c(1)}\leq\bar{Q}_l^t+b_l v_{\mathrm{up}(l)},\forall\, l\in\mathcal{L}.$$

These show Equations (27) and (28).

**2- Induction step**

Assume the statements in *Lemma III* are true for some $n\geq 1$. We now show it also holds for $n+1$.

**a)** Consider first some fixed $l\in\mathcal{L}$. Define $\Phi_l$ by $f_l\triangleq\phi_l\left(P_l^t,Q_l^c,v_{\mathrm{up}(l)}\right)$ from Equation (8.c). We have

$$\mathrm{grad}(\phi_l)=\begin{pmatrix}\dfrac{2P_l^t}{v_{\mathrm{up}(l)}}\\ \dfrac{2(Q_l^c)}{v_{\mathrm{up}(l)}}\\ -\dfrac{(P_l^t)^2+(Q_l^c)^2}{\left(v_{\mathrm{up}(l)}\right)^2}\end{pmatrix} \quad (61).$$

Define $M(t)$ for $t\in[0,1]$ as

$$M(t)=t\begin{pmatrix}P_l^{t(n)}\\ Q_l^{c(n)}\\ v_{\mathrm{up}(l)}^{(n)}\end{pmatrix}+(1-t)\begin{pmatrix}P_l^{t(n-1)}\\ Q_l^{c(n-1)}\\ v_{\mathrm{up}(l)}^{(n-1)}\end{pmatrix}.$$

Then by Equation (8.c) and by the fundamental law of calculus

$$f_l^{(n+1)}-f_l^{(n)}=\Phi_l\big(M(1)\big)-\Phi_l\big(M(0)\big)=\begin{pmatrix}\Delta P_l^{t(n)}\\ \Delta Q_l^{c(n)}\\ \Delta v_{\mathrm{up}(l)}^{(n)}\end{pmatrix}.\int_0^1\mathrm{grad}\Phi_l\big(M(t)\big)dt.$$

We first bound each component of the gradient. For $0\le t\le 1$, by the induction property (23), (27), and (28) at $n-1$ and $n$ (note that $\omega$ is a feasible solution thus $P_l^t\le \bar{P}_l^t\le P_l^{\max}$, $Q_l^c=Q_l^t+b_l v_{\mathrm{up}(l)}\le \bar{Q}_l^t+b_l v^{\max}\le Q_l^{\max}+b_l v^{\max}$

$$(1-t)P_l^{t(n)}+tP_l^{t(n-1)}\le P_l^{\max}$$
$$(1-t)Q_l^{c(n)}+tQ_l^{c(n-1)}\le Q_l^{\max}+b_l v^{\max}$$
$$v_{min}\le (1-t)v_l^{(n)}+tv_l^{(n-1)}.$$

Furthermore, $f^n\ge 0$ for any $n\ge 0$ and the matrices in Equations (47) and (35) are non-negative, therefore, for all $n\ge 0$:

$$\hat{P}_l^t\le P_l^{t(n)}$$
$$\hat{Q}_l^c\le Q_l^{c(n)}$$

and thus

$$(1-t)P_l^{t(n)}+tP_l^{t(n-1)}\ge \hat{P}_l^t$$
$$(1-t)Q_l^{c(n)}+tQ_l^{c(n-1)}\ge \hat{Q}_l^c.$$

By Equation (61) it follows that (entry-wise):

$$\left|\mathrm{grad}\big(\Phi_l(M(t)\big)\right|\le\begin{pmatrix}2\pi_l\\ 2\varrho_l\\ \vartheta_l\end{pmatrix}$$

thus, we have

$$\left|\Delta f_l^{(n+1)}\right| \leq 2\pi_l \left|\Delta P_l^{t(n)}\right| + 2\varrho_l \left|\Delta Q_l^{c(n)}\right| + \vartheta_l \left|\Delta v_{\mathrm{up}(l)}^{(n)}\right|.$$

Now this is true for some $l$, so in matrix form we have:

$$\left|\Delta f^{(n+1)}\right| \leq 2\mathrm{diag}(\pi)\left|\Delta P^{(n)}\right| + 2\mathrm{diag}(\varrho)\left|\Delta Q^{c(n)}\right| + \mathrm{diag}(\vartheta)\left|\Delta v^{(n)}\right| \qquad (62).$$

By the construction of $P^{(n)}$, $Q^{c(n)}$, and $v^{(n)}$ we have

$$\left|\Delta P^{(n)}\right| \leq \mathbf{H}\mathrm{diag(r)}\left|\Delta f^{(n)}\right| \qquad (63)$$

$$\left|\Delta Q^{c(n)}\right| \leq \mathbf{F}\left|\Delta f^{n}\right| \qquad (64)$$

$$\left|\Delta v^{(n)}\right| \leq \mathbf{D}\left|\Delta f^{(n)}\right| \qquad (65)$$

combined with (7). and (62) this gives

$$\left|\Delta f^{(n+1)}\right| \leq \mathbf{E}\left|\Delta f^{(n)}\right| \qquad (66)$$

applying the induction property (21) it comes

$$\left|\Delta f^{(n+1)}\right| \leq \mathbf{E}^{n}\left|\Delta f^{(1)}\right| \qquad (67)$$

which shows that he induction property (21) also holds for $n + 1$.

■

**b)** We have

$$\Delta v^{(1)} = (-\mathbf{D})\Delta f^{(1)} \qquad (68)$$

and we have already noted that $\Delta f^{1} \leq 0$, thus

$$\left|\Delta v^{(1)}\right| = \Delta v^{(1)} = (\mathbf{D})\left|\Delta f^{(1)}\right|.$$

Using (43) we have

$$\left|\Delta v^{(n+1)}\right| \leq \mathbf{D}\left|\Delta f^{(\mathrm{n}+1)}\right|.$$

apply (65), (67) and *C3*

$$\left|\Delta v^{(n+1)}\right| \leq \eta^n \mathbf{D}\left|\Delta f^{(1)}\right|. \tag{69}$$

apply (68)

$$\left|\Delta v^{(n+1)}\right| \leq \eta^n\left|\Delta v^{(1)}\right| \tag{70}$$

it follows that (22) also holds for $n+1$.

■

Furthermore we have

$$\left|v^{(n+1)}-v^{(1)}\right| \leq \left|\Delta v^{(n+1)}\right| + \cdots + \left|\Delta v^{(2)}\right| \leq (\eta^n + \cdots + \eta)\left|\Delta v^{(1)}\right| \leq \frac{\eta}{1-\eta}\left|\Delta v^{(1)}\right| \leq \left|v^{(1)}-v\right|$$

thus (noting that $v^{(1)} \geq v$)

$$v^{(1)}-v^{(n+1)} \leq \left|v^{(n+1)}-v^{(1)}\right| \leq \left|v^{(1)}-v\right| = v^{(1)}-v$$

thus

$$v^{(n+1)} \geq v \tag{71}$$

It follows that (23) also holds for $n+1$.

■

**c)** We have

$$\Delta P^{(1)} = \mathbf{H}\mathrm{diag}(r)\Delta f^{(1)} \tag{72}$$

and we have already noted that $M$ entries of $\Delta f^{(1)}$ are non-zero ($\Delta f_m^{(1)} < 0, \forall\, m \in \mathcal{M}$), thus:

$$\left|\Delta P_l^{t(1)}\right| = -\Delta P_l^{t(1)} = -\sum_{m\in\mathcal{M}}\left(\left(\mathbf{H}\mathrm{diag}(r)\right)_{lm}\Delta f_m^{(1)}\right) = \sum_{m\in\mathcal{M}}\left(\left(\mathbf{H}\mathrm{diag}(r)\right)_{lm}\left|\Delta f_m^{(1)}\right|\right) \tag{73}.$$

Using (47) and (67) we have

$$\left|\Delta P_l^{t(n+1)}\right| \leq (\mathbf{H}\mathrm{diag}(r))_l\left|\Delta f^{(n+1)}\right| \leq \left(\mathbf{H}\mathrm{diag}(r)\mathbf{E}^n\left|\Delta f^{(1)}\right|\right)_l = \sum_{m\in\mathcal{M}}\left((\mathbf{H}\mathrm{diag}(r)\mathbf{E}^n)_{lm}\left|\Delta f_m^{(1)}\right|\right) \tag{74}.$$

Applying *C5* for $n$ times we have

$$\mathbf{H}\text{diag}(r)\ \mathbf{E}^n \leq (\eta)^{n-1}\mathbf{H}\text{diag}(r)\mathbf{E}, \forall\ n \geq 1 \tag{75}$$

thus

$$(\mathbf{H}\text{diag}(r)\ \mathbf{E}^n)_{lm} \leq (\eta)^{n-1}(\mathbf{H}\text{diag}(r)\mathbf{E})_{lm}, \forall\ n \geq 1, \quad \forall\ m \in \mathcal{M} \tag{76}$$

Furthermore from *C4* for $l \in \mathcal{L}^{\mathcal{M}}$ we have

$$(\mathbf{H}\text{diag}(r)\ \mathbf{E})_{lm} \leq \eta\big(\mathbf{H}\text{diag}(r)\big)_{lm} \quad \forall\ m \in \mathcal{M} \tag{77}$$

combining with (76)

$$(\mathbf{H}\text{diag}(r)\ \mathbf{E}^n)_{lm} \leq (\eta)^n\big(\mathbf{H}\text{diag}(r)\big)_{lm}, \forall\ l \in \mathcal{L}^{\mathcal{M}}, \forall\ m \in \mathcal{M} \tag{78}$$

combining with (74)

$$\begin{aligned}\left|\Delta P_l^{t(n+1)}\right| &\leq \big(\mathbf{H}\text{diag}(r)\big)_l |\Delta f^{n+1}| \leq \sum_{m\in\mathcal{M}} \left((\mathbf{H}\text{diag}(r)\mathbf{E}^n)_{lm} \left|\Delta f_m^{(1)}\right|\right) \\ &\leq \eta^n \sum_{m\in\mathcal{M}} \left(\big(\mathbf{H}\text{diag}(r)\big)_{lm} \left|\Delta f_m^{(1)}\right|\right), \qquad \forall\ l \in \mathcal{L}^{\mathcal{M}}\end{aligned} \tag{79}$$

apply (73)

$$\left|\Delta P_l^{t(n+1)}\right| \leq \eta^n \left|\Delta P_l^{t(1)}\right|, \qquad \forall\ l \in \mathcal{L}^{\mathcal{M}} \tag{80}$$

it follows that (24) also holds for $n + 1$.

■

Furthermore for $l \in \mathcal{L}^{\mathcal{M}}$ we have (noting that $\eta < 0.5$)

$$\begin{aligned}\left|P_l^{t(n+1)} - P_l^{t(1)}\right| &\leq \left|\Delta P_l^{t(n+1)}\right| + \cdots + \left|\Delta P_l^{t(2)}\right| \leq (\eta^n + \cdots + \eta) \left|\Delta P_1^{t(1)}\right| \leq \frac{\eta}{1-\eta}\left|\Delta P_l^{t(1)}\right| \\ &\leq \left|P_l^{t(1)} - P_l^t\right|, \\ &\qquad \forall\ l \in \mathcal{L}^{\mathcal{M}}\end{aligned} \tag{81}$$

thus (noting that $P^{(1)} \leq P$)

$$P_l^{t(n+1)} - P_l^{t(1)} \leq \left|P_l^{t(n+1)} - P_l^{t(1)}\right| \leq \left|P_l^{t(1)} - P_l^{t}\right| = P_l^{t} - P_l^{t(1)}, \forall\, l \in \mathcal{L}^{\mathcal{M}} \tag{82}$$

thus $P_l^{t(n+1)} \leq P_l^t\ , \forall\, l \in \mathcal{L}^{\mathcal{M}}$.

It follows that (25) also holds for $n+1$.

■

**d)** We have

$$\hat{P}_l^t \leq P_l^{t(n)} \leq \bar{P}_l^t\ , \forall\, l \in \mathcal{L} \tag{83}$$

$$\hat{Q}_l^c \leq Q_l^{c(n)} = Q_l^{t(n)} + b_l v_{\mathrm{up}(l)}^{(n)} \leq \bar{Q}_l^t + b_l v_{\mathrm{up}(l)}\ , \quad \forall\, l \in \mathcal{L}$$

thus from step 1 of *Algorithm I,* Equation (11.e) and noting that $v^{(n)} \geq v$, we have

$$f_l^{(n+1)} \leq \bar{f}_l\ , \forall\, l \in \mathcal{L} \tag{84}$$

it follows that (26) also holds for $n+1$.

■

**e)** We have already shown that $f_l^{(n+1)} \leq \bar{f}_l\ , \forall\, l \in \mathcal{L}$ and $v^{(n+1)} \geq v$ in Equations (84) and (71) respectively. From Equations (36), (37), and (52) one can show that

$$P_l^{t(n+1)} \leq \bar{P}_l^t\ , \forall\, l \in \mathcal{L} \tag{85}$$

$$Q_l^{c(n+1)} \leq \bar{Q}_l^t + b_l v_{\mathrm{up}(l)}\ , \forall\, l \in \mathcal{L} \tag{86}$$

It follows that (27) and (28) also hold for $n+1$.

Both basis and inductive steps are proved which completes the proof of *Lemma III.*

■